\documentclass[mathpazo]{cicp}

\usepackage{amssymb}
\usepackage{amsmath}
\usepackage{amsthm}
\usepackage{braket}
 
\usepackage{graphicx}
\usepackage[style=base]{caption}
\usepackage{subcaption}

\usepackage{comment}
\usepackage{listings}
\usepackage{algorithm2e}
\usepackage[usenames,dvipsnames]{xcolor}

\newcommand{\bm}[1]{\mathbf{#1}}
\newcommand{\all}{\forall}
\usepackage{hyperref}

\newcommand{\average}[1]{\ensuremath{\langle#1\rangle} }
\newcommand{\jump}[1]{\ensuremath{[\![#1]\!]} }

\begin{document}
\title{Capturing near-equilibrium solutions: a comparison between high-order discontinuous Galerkin methods and well-balanced schemes}

 \author[Han Veiga et.~al.]{Maria Han Veiga\affil{1}\comma \affil{2}\comma\corrauth, David A. Romero Velasco\affil{1}\comma\affil{3}\comma\affil{4}
       R\'emi Abgrall\affil{2}, and Romain Teyssier\affil{1}}
 \address{\affilnum{1}\ Institute of Computational Science,
          University of Zurich, Switzerland \\
           \affilnum{2}\ Institute of Mathematics,
           University of Zurich, Switzerland\\
           \affilnum{3}\ Universidad Aut\'{o}noma del Estado de Morelos, Mexico\\
           \affilnum{4}\ Instituto de Ciencias F\'{i}sicas, Universidad Nacional Aut\'{o}noma de M\'{e}xico, Mexico}
 \emails{{\tt mariaioque.hanveiga@math.uzh.ch} (M.~Han Veiga), {\tt david.velasco@icf.unam.mx} (D.~Velasco), {\tt remi.abgrall@math.uzh.ch} (R.~Abgrall),
          {\tt romain.teyssier@uzh.ch} (R.~Teyssier)}


\begin{abstract}
Equilibrium or stationary solutions usually proceed through the exact balance between hyperbolic transport terms and source terms. 
Such equilibrium solutions are affected by truncation errors that prevent any classical numerical scheme from capturing the evolution of small amplitude waves of physical significance. 
In order to overcome this problem, we compare two commonly adopted strategies: going to very high order and reduce drastically the truncation errors on the
equilibrium solution, or design a specific scheme that preserves by construction the equilibrium exactly, the so-called well-balanced approach. 
We present a modern numerical implementation of these two strategies and compare them in details, using hydrostatic but also dynamical equilibrium solutions of several simple test cases. 
Finally, we apply our methodology to the simulation of a protoplanetary disc in centrifugal equilibrium around its star and model its interaction with an embedded planet, illustrating in a realistic application the strength of both methods. 
\end{abstract}


\maketitle


\section{Introduction}
\label{sec1}
Hyperbolic balance laws are used to describe many dynamical problems in natural sciences. 
They are defined as a set of conservation laws with associated source terms, 
which model the production or destruction of the corresponding conserved quantity. 
Many physical systems of scientific interest can be described by a system of hyperbolic conservation laws with source terms, or in short, hyperbolic balance laws. 

Hyperbolic balance laws are particularly challenging because they feature equilibrium solutions that result from the exact cancellation of the left hand side and the right hand side of these equations.
Small truncation errors can perturb this equilibrium solution, leading to the production of spurious waves that can dominate over
the real waves  that control the physics of the problem at hand.

For example, for the case of the inviscid Euler equations with a gravity source term (also known as the Euler-Poisson system), hydrostatic steady states are important in, for example, hydraulics \cite{pares2007, berthon2012, bermudez2017} and astrophysics \cite{supernova1996, solarbook, galaxybook}. The difficulty here is to capture properly sound waves, gravity waves or convective flows, whose amplitude can be comparable to the truncation errors of a second order method and a reasonable grid resolution. 

General steady states with non constant velocity fields are also found to be important in planetary sciences, namely in the early stages of protoplanetary discs, where the source term models the gravity of a central star\cite{surville2016}, and are balanced by the centrifugal and pressure forces. The challenge here is to be able to resolve the interaction of a small planet with the gaseous disc, leading to the formation of a small amplitude spiral wave that can be dominated by the truncation errors of the equilibrium solution. In this context, the classical approach is to use a cylindrical mesh, reducing drastically discretisation errors along circular orbits. It is however desirable to find a solution on a Cartesian mesh, as it allows to deal with more general cases which are not strictly axisymmetric. 

In summary, solving for such flows which are close to equilibrium can be very challenging for a naive, low order numerical method on a mesh not necessarily adapted to the geometry of the equilibrium solution as the truncation error incurred while solving the steady state can be larger than the small amplitude waves of interest. 

There are nowadays many practical numerical methods with very low truncation errors. A class of such methods are the so-called discontinuous Galerkin (DG) methods \cite{ern}. These methods, at least for smooth and regular problems, can be made as accurate as desired. This means that, at least in principle, the amplitude of the truncation errors can be reduced to an arbitrarily small value.
This requires an appropriate way to implement the source terms in the DG formalism \cite{praveen, gangli}. This also requires the use of a high enough resolution mesh to capture the equilibrium solution, which translates into higher computational cost for higher order solutions. 

There is another strategy that allows one to use a low-order method, while capturing almost exactly the equilibrium solution. This is called the well-balanced approach (introduced in detail \cite{leveque2006}), which is concerned with numerical schemes that satisfy the discrete equivalent of an underlying steady state, effectively, taking into account the existence of a steady state (or near steady state) solution. 

The natural question is thus whether exact well-balancedness is required in practice or if methods that solve the PDE (including the source term) need only to be very accurate. This is the question we wish to explore in this paper on several examples of interest for natural sciences in general and astrophysics in particular.

Let $d$, $e \in \mathbb{N}$, $\Omega$ an open subset of $\mathbb{R}^e$ and $\bm{f}_j$ for $1\leq j \leq d$ be smooth functions from from $\Omega$ into $\mathbb{R}^e$. A general $e$-size system of $d$-dimensional hyperbolic balance laws can be written in the following form: 

\begin{equation}
\label{eq:main}
    \frac{\partial \bm{w}}{\partial t} + \sum_{j=1}^{d} \frac{\partial}{\partial x_j} \bm{f}_j(\bm{w}) - \bm{s}(\bm{w},\bm{x}) = 0, \quad \bm{x}=(x_1,...,x_d)\in \mathbb{R}^d, \quad t > 0,
\end{equation}

where the vector valued function $\bm{w}=(w_1,...,w_e):\mathbb{R}^d \times [0,\infty)\to \Omega$ denotes the solution, the functions $\bm{f}_j=[f_{1j},...,f_{ej}]^T$ are flux-functions and $\bm{s}(\bm{w},\bm{x})$ is the vector of source terms. We denote vectors in bold $\bm{v}$ and a component of the vector as $v$, where the index is omitted if not important. 

In order to solve hyperbolic balance laws, one can use classical methods for hyperbolic conservation laws  (i.e. when $\textbf{s}(\textbf{w},\textbf{x})=\textbf{0}$) in conjunction with an operator-split approach to add the source terms. However, problems can arise when one tries to model flows near equilibrium states, for which (\ref{eq:main}) admits a steady state solution such that:

\begin{equation}
\label{eq:fluxsourcebal}
\sum_{j=1}^{d} \frac{\partial}{\partial x_j} \textbf{f}_j(\bf{w}) - \bf{s}(\bf{w},\bf{x}) = 0
\end{equation}

In this work, we are mainly interested in solving the Euler-Poisson system with an analytical gravitational potential $\Phi$ with moving steady states (where the velocity field $\bm{v}\not\equiv \bm{0}$). We restrict ourselves to one- or two-dimensional cases given by $ e = 3$ and $d = 1$, or $e = 4$ and  $d = 2$ respectively. We however present the main equations in the 2D case only, as given in \eqref{eq:euler}. 

\begin{align}
\label{eq:euler}
    \frac{\partial \bm{w}}{\partial t} + \sum_{j=1}^{2} \frac{\partial}{\partial x_j} \bm{f}_j(\bm{w}) - s(\bm{w},\bm{x}) = 0
\end{align}
where

\begin{align*}
 \textbf{w}=\begin{bmatrix}
         \rho \\
         \rho v_x \\
         \rho v_y \\
         E
        \end{bmatrix}, \, \textbf{f}_1(\textbf{w}) = \begin{bmatrix}
         \rho v_x  \\
         \rho v_x^2 + p  \\
         \rho v_x v_y \\
         v_x(E+p) 
        \end{bmatrix}, \,
        \textbf{f}_2(\textbf{w}) = \begin{bmatrix}
         \rho v_y  \\
         \rho v_x v_y \\
         \rho v_y^2 + p \\
         v_y(E+p)
        \end{bmatrix}, \,
        \textbf{s}(\textbf{w}) = \begin{bmatrix}
         0 \\
         -\rho\frac{\partial}{\partial x} \Phi \\
         -\rho\frac{\partial}{\partial y} \Phi \\
         -\rho \bm{v}\cdot \nabla \Phi
        \end{bmatrix}.
\end{align*}

Here $\rho$ is the mass density, $\bm{v}=(v_x,v_y)$ the velocity and $E$ the total energy given by the sum of internal and kinetic energy.

\[E= \rho\epsilon + \frac{1}{2} \rho |\bm{v}|^2.\] 

In addition, $p$ denotes the pressure and we assume $p=p(\rho,e)$ is a known function. Furthermore, we assume an ideal gas, such that the system is closed with the equation of state:

\[p = \rho \epsilon(\gamma - 1), \]

where $\gamma$ denotes the adiabatic index. The source terms, shown on the right hand side of the momentum and energy equations, model the effect of the gravitational forces on the fluid, for a given potential $\Phi$.

To correctly solve these equations numerically and capture small perturbations to the steady state, dedicated computational methods are required to solve the discrete version of the source-flux balance \eqref{eq:fluxsourcebal}. For non well-balanced methods, there is no guarantee that the truncation errors induced by discretising  the steady state solution are not greater than the small perturbations we want to describe.

The design of \emph{well-balanced schemes} (i.e. schemes which satisfy exactly a discrete equivalent of the underlying steady state) has been an active field of research, first coined in \cite{Leroux96}. There have been many attempts to deal with this aspect, in particular for the shallow water equations, where steady states can represent the lake at rest case (hydrostatic equilibrium) \cite{bermudez94} or a running river (non-trivial velocity equilibrium state)\cite{Ricchiuto15}. 

For the Euler-Poisson system  there have been several recent contributions. We do not intend to give an exhaustive account of all the work that has be done in this topic, but we refer to: \cite{Kapelli} where the authors design a well-balanced first and second order accurate finite volume scheme for approximating the Euler equations with gravitation  using a discretisation of the hydrostatic equilibrium for the pressure reconstruction, \cite{gangli,gangli2} where a similar approach to treat hydrostatic equilibria, isothermal and polytropic equations of state achieves a high order well-balanced discontinuous Galerkin scheme, and \cite{Klingenberg} where a relaxation scheme is adopted.  

In \cite{surville2016} a second order finite volume method dealing with non zero velocities is presented in the context of protoplanetary discs. Concerning more general classes of steady states, the survey \cite{noelle2016} describes two classes of schemes, one based on high-order accurate, non-oscillatory finite difference operators which are well-balanced for a general class of equilibria, and another one based on well-balanced quadratures, showing the suitability of these methods on the Shallow Water equations, and the work in \cite{Pares2006}, describing a high order path-conservative scheme and well balanced reconstruction, however, the analysis for this work is restricted to 1-dimension quasi-linear hyperbolic systems.

On the other hand, due to the tractability of modern, very-high-order methods, one could ask whether these methods alone could be enough to solve the equilibrium solutions to a high enough accuracy.

In this paper, we make a comparative study between a new method which is truly well balanced and a popular class of very-high-order methods, namely the discontinuous Galerkin method. We would like to answer the following fundamental questions, considering both hydrostatic and moving equilibria solutions:
\begin{enumerate}
	\item  Are there cases where using a high order scheme is sufficient to capture solutions close to a steady state?
	\item  Under which circumstances is it necessary to use a well balanced method?
	\item What is the cost associated to each approach and how does it balance with accuracy?
\end{enumerate}

In particular, we compare a well-balanced, high-order discontinuous Galerkin method with a non well-balanced, high-order discontinuous Galerkin method under different steady state regimes, both hydrostatic and stationary (with a non-zero velocity). 

The outline of this paper is as follows: a brief introduction on equilibrium solutions, as well as the description of the Runge-Kutta discontinuous Galerkin (RKDG) method is provided in section 2. In section 3, we describe our well balanced formulation of RKDG. In section 4, a set of benchmark problems are defined, both in one and two dimensions, and quantitative results are presented, followed by our final discussion in section 5.

\section{Preliminaries}
\label{sec1}
\subsection{Steady state solutions}
A solution  $\bm{w}$ is said to be a steady state solution of (\ref{eq:main}) if it fulfils the following relation
\begin{equation}
    \label{eq:generalfluxsource}
\sum_{j=1}^{d} \frac{\partial}{\partial x_j} \textbf{f}_j(\bf{w}) - \bf{s}(\bf{w},\bf{x}) = 0.
\end{equation}
for $w:\mathbb{R}^d \times  [0, \infty ) \to \Omega$,  $\textbf{f}_j:\Omega \to \mathbb{R}^e$ and $s:\Omega\times\mathbb{R}^d \to \mathbb{R}^e$, $\Omega \subset \mathbb{R}^e$.

 We call $\bm{w}$  a hydrostatic steady state of (\ref{eq:euler}) if the pressure component fulfils the following relation
\begin{equation}
    \nabla p = -\rho\nabla \Phi,
\end{equation}
for $\rho, p:\mathbb{R}^d \times  [0, \infty ) \to \Omega$, and $\Phi:\Omega\times\mathbb{R}^d \to \mathbb{R}^e$, $\Omega \subset \mathbb{R}^e$ and a gravity potential $\Phi\in\mathcal{C}^1$.

 Let the following be the standardised form of an time explicit numerical scheme for \eqref{eq:euler}, where $H(\cdot)$ denotes the update function for each timestep $n$ in a quantity indexed by $k$ (e.g. cell average) and $q,~p$ denote the stencil size:
\begin{equation}
\label{eq:genericupdate}
    w^{n+1}_k = w^n_k + \frac{\Delta t}{\Delta x} H(w^n_{k-q}. ..., w^n_{k+p})
\end{equation}

The numerical scheme is exactly well-balanced if for a steady state solution $w$, the following holds
\begin{equation}
    H(w^n_{k-q}, ..., w^n_{k+p}) = 0.
\end{equation}

The scheme is said to be well-balanced with order $N_p$ if, for a steady state solution $w$, the following holds
\begin{equation}
    |H(w^n_{k-q}, ..., w^n_{k+p})| = \mathcal{O}(\Delta x^{N_p+1}).
\end{equation}

A formal definition of the above is given in \cite{Pares2006}.

\theoremstyle{remark}
\begin{remark} 
In the astrophysics literature, in particular in the planet formation community \cite{armitageplanet}, the following equilibrium relation between the centrifugal force associated to the cross-radial velocity $v_{\theta}$, the gradients of the thermal pressure $p$ and the gradient of the gravitational potential $\Phi$ is often referred to as a \emph{dynamical equilibrium}:
\begin{equation}
    \frac{v_{\theta}^2}{r} = \frac{1}{\rho}\nabla p + \nabla \Phi
\end{equation}
It is often the case that the pressure gradient $\nabla p$ is assumed to be small and thus can be neglected \cite{armitageplanet}. For the purpose of this work, we considered initial conditions which are strictly hydrostatic, and also more general \emph{stationary solutions} or \emph{steady states} of the Euler-Poisson equation.
\end{remark}

\subsection{Runge Kutta Discontinuous Galerkin (RKDG) method}
\label{subsec1}

Consider a regular domain $D \in \mathbb{R}$, approximated by $K$ non-overlapping elements such that $\bigcup_{K\in D_h} K \approx D$. The 2-dimensional tessellation is given by the tensor product of the 1-dimensional discretizations, thus yielding square volumes (or cubic volumes in 3D). Let $T_h$ denote the Cartesian tessellation of the domain $D$ where our problem is defined.

We seek for the approximate solution $w_h(t)$ in the finite element space of discontinuous functions $V_h$:
\[ V_h = \{v_h \in \mathcal{L}^{\infty}(D): v_h|K \in V_h(K), \all K \in T_h \}. \]
We take $V_h(K)$ to be the collection of polynomials of at most degree $N_p$. 

Following the Runge Kutta discontinuous Galerkin (RKDG) method described in \cite{cs5}, we write the weak formulation for each component of \eqref{eq:main} by multiplying the system by a smooth test function $v(x)$ and integrate over a control volume $K$:

\begin{align}
\label{eq:weakform}
\begin{split}
\frac{\mbox{d}}{\mbox{dt}}\int_K w(\bm{x},t)v(\bm{x})\mbox{d\bf{x}} + \sum_{e\in\partial K}\int_e f(w(\bm{x},t))\cdot n_{e,K} v(\bm{x}) \mbox{d}\Gamma &- \int_K f(w(\bm{x},t))\cdot \nabla v(\bm{x})\mbox{d\bf{x}}\\ &= \int_K s(w(\bm{x},t)) v(\bm{x})\mbox{d\bf{x}}
\end{split}
\end{align}
for any smooth $v(\bm{x})$. We denote the outward unit normal as $n_{e,K}$ and edge as $e$.

The following integrals are approximated with a suitable order numerical quadrature (where $\{\bm{x}_i,\omega_i\}_{i=0}^{M,L}$ denotes the set of quadrature points and weights): 

\begin{align}
 &\int_e f(w(\bm{x},t))\cdot n_{e,K} v(\bm{x}) \mbox{d}\Gamma   \approx \sum_{i=0}^L f(w(\bm{x}_i,t))\cdot n_{e,K} v(\bm{x}_i) \omega_i |e| ,  \\
 &\int_K f(w(\bm{x},t))\cdot \nabla v(\bm{x})\mbox{d\bf{x}}  \approx \sum_{j=0}^M f(w(\bm{x}_j,t))\cdot \nabla v(\bm{x}_j) \omega_j |K| , \\
 &\int_K s(w(\bm{x},t)) v(\bm{x})\mbox{d\bf{x}}  \approx \sum_{j=0}^M s(w(\bm{x}_j,t)) v(\bm{x}_j) \omega_j |K| .
\end{align}

Note that there is an ambiguity in the definition of the flux $f(w(\bm{x},t))\cdot n_{e,K}$  since $w$ can be multi-valued at the cell interface. In order to overcome this inconsistency, this term is replaced by a single-valued numerical flux $h_{e,K}(\bm{x},t)$ computed using a Riemann solver. The exact solution is replaced by the finite dimensional approximate solution $w_h = \sum_{i=0} \hat{w}_i(t)\psi_i (\bm{x}) $, where $\hat{w}_i(t)$ is given by the $L^2$ inner product between $w(\bm{x},t)$ and $\psi_i$, a basis element of $V_h (K)$, and the test functions $v(\bm{x})$ are replaced by $v_h(\bm{x}) \in V_h(K)$. This yields the following numerical scheme:

\begin{equation*}
\begin{split}
 w_h(t=0)&=P_{V_h}(w_0) \\
 \frac{\mbox{d}}{\mbox{dt}}\int_K w_h(\bm{x},t)v_h(x)\mbox{d\bf{x}} &= -\sum_{e\in\partial K}\sum_{i=1}^L h_{e,K}(\bm{x}_i,t)v_h(\bm{x}_i) \omega_i |e| \nonumber \\ &
 \quad + \sum_{j=1}^M f(w_h(\bm{x}_j,t))\cdot \nabla v_h(\bm{x}_j) \omega_j |K| \\
 &\quad +\sum_{j=1}^M S(w_h(\bm{x}_j,t)) v_h(\bm{x}_j) \omega_j |K| \qquad \forall v_h(\bm{x})\in V(K), \forall K \in T_h, \label{eq:dg}
\end{split}
\end{equation*}
where operator $P_{V_h}$ denotes the $L^2$ projection of the initial data $w_0(x)$ into the space of finite elements $V_h$. 

In addition, throughout this work, we make the following choices:
\begin{enumerate}
    \item{ We denote by $\{ \psi \}_{i=0}^{N_p}$ the Legendre basis vectors spanning $V_h(K)$, subject to the following normalisation:
    \[ \int_{-1}^{1} \psi_i(x)\psi_j(x) dx = \delta_{ij}\, ; \] }
    \item{ We take the numerical quadrature points and weights $\{ x_i, \omega_i \}_{i=0}^M $ to be Gauss-Legendre quadrature points;}
    \item{ We use local Lax-Friedrichs flux as the numerical flux. We note that our analysis works for any consistent numerical flux function that is Lipschitz continuous in both arguments, non-decreasing in its first argument and non-increasing in its second argument.}
\end{enumerate}

\subsubsection{Time discretisation}
We use the TVD Runge-Kutta time discretization as in \cite{gottlieb}. Let $\{t^n\}^N_{n=0}$ be a partition of $[0,T]$ and $\Delta t^n = t^{n+1}-t^{n}$, $n=0, ... , N-1$, then the time marching algorithm is given in algorithm \ref{algo:time}. 

\begin{algorithm}
\SetAlgoLined
\KwData{$w^0_h = P_{V_h}(w_0)$}
\KwResult{$w^{n+1}_h$}
\For{n = 0, ... N-1}{
    $w^{(0)}_h = w^n_h$\\
    \For{ i = 1, ... k+1}{
        $k_i = \mathcal{L}(t^n + c_i, w_h^{(0)} + h \sum_{j=1}^{i-1} a_{i,j} k_j)$
    }
    $w^{n+1}_h = w_h^{(0)} + h\sum_{i=1}^s b_i k_i$
}

 \caption{TVD RK time marching algorithm}
 \label{algo:time}
\end{algorithm}

The parameters $a_{i,j}$, $b_i$ and $c_{i}$ can be found in Tables \ref{butcher}, see \cite{gottlieb}.

Given that an explicit time integrator is used, the timestep $\Delta t$ has to fulfill a Courant-Friedrich-Lewy (CFL) condition to achieve numerical stability. In this work, the timestep $\Delta t^K$ at cell $K$ is calculated as \cite{shu89}. Furthermore, the introduction of a source term can introduce additional constraints on the timestep. As described in \cite{zhang}, the timestep for a solution approximation of degree at most $N_p$, we choose the minimum of the expression below: 

\[
\Delta t^K = \min \left( \frac{C}{2N_p + 1} \left ( \sum_{i=1}^d \frac{|v_i^K| + c_s^K}{\Delta x_i^K} \right)^{-1}, \frac{1}{\sqrt{2\gamma(\gamma - 1)}} \frac{c_s^K}{|\nabla \Phi^K|} \right)
\]

where $c_s = \sqrt{\gamma p / \rho}$ is the sound speed, $v_i^K$ is the $i^{th}$ component of the velocity average at cell $K$, $\Delta x_i^K$ is the mesh-width in the $i^{th}$ dimension and $|\nabla \Phi|^K$ the magnitude of the gradient of $\Phi$ at cell $K$. The constant $C$ is chosen to be small, for example, \textbf{0.2}.

\begin{table}[h]
\begin{center}
\begin{tabular}{cc}
$
\begin{array}{c|cc}
0\\
1/2 & 1/2\\
\hline
& 1/2 & 1/2 \\
\end{array}
$
& 
$\begin{array}{c|ccc}
0\\
1 & 1 & & \\
3/4& 1/4 & 1/4\\
\hline
&1/6 & 1/6 & 1/3\\
\end{array}$
\\ &\\
SSP(2,2)  &  SSP(3,3)
\end{tabular}
\end{center}
\caption{\label{butcher} Runga-Kutta Butcher tableaus for the TVDRK schemes.}
\end{table}

\begin{table}[h]
\begin{center}
\scalebox{0.7}{
\begin{tabular}{c}
$\begin{array}{c|ccccc}
0\\
0.39175222700392 & 0.39175222700392 & & & & \\
0.58607968896779 & 0.21766909633821  & 0.36841059262959 & & & \\
0.47454236302687 & 0.08269208670950  & 0.13995850206999 &  0.25189177424738 &  & \\
0.93501063100924 & 0.06796628370320  & 0.11503469844438 &  0.20703489864929 & 0.54497475021237 & \\
\hline
& 0.14681187618661 & 0.24848290924556 & 0.10425883036650 & 0.27443890091960 &  0.22600748319395\\
\end{array}$
\\ 
SSP(4,5)
\end{tabular}
}
\end{center}
\caption{\label{butcher} Runga-Kutta Butcher tableaus for the TVDRK schemes.}
\end{table}

\subsection{Discretisation error of DG for Steady Euler System}

In this short section we show that the traditional Runge Kutta discontinuous Galerkin method is inherently not well balanced and specify the source of approximation error for each conserved variable. 

We consider the 1-dimensional Euler system \eqref{eq:11deuler} and the DG discretisation described above. 

\begin{subequations}\label{eq:11deuler}
     \begin{align}
      &\frac{\partial}{\partial t} \rho +  \frac{\partial}{\partial x} \rho v = 0 \label{eq:1ddensity}\\
      &\frac{\partial}{\partial t} \rho v +  \frac{\partial}{\partial x}( \rho v^2 + p) = -\rho \frac{\partial}{\partial x} \Phi \label{eq:1dmomentum}\\
      &\frac{\partial}{\partial t} E +  \frac{\partial}{\partial x}(v(E+p)) = -\rho v \frac{\partial}{\partial x} \Phi  \label{eq:1denergy}
     \end{align}
\end{subequations}
    
Considering, for example, \eqref{eq:1ddensity} in some control volume $K$ mapped to $[-1,1]$ interval, and modal coefficient $i$:
\begin{equation}
    \partial_t \tilde{\rho_i} = \int_K \rho v \partial_x \psi_i \mbox{dx} - \oint_{\partial K} \rho v\psi_i \cdot n \mbox{d}\Gamma
\end{equation}

If the solution of the system is a steady state, then $\partial _t \tilde{\rho_i} = 0 \quad \forall i$. 

One can write the update $H_\rho$ vector (as in \eqref{eq:genericupdate}), where each component corresponds to $i^{th}$ the modal update, with the associated test function $\psi_i$:
\[H_\rho^{i} = \sum_j \rho v \partial_x \psi_i(x_j)w_j - \hat{\rho v}\psi_i(1) + \hat{\rho v}\psi_i(-1)\]
To evaluate $\hat{\cdot}$ we need an approximate flux function to combine the left and right hand-side values of the flux. We consider the Lax-Friedrichs flux, for example, defined as:
\[ \hat{f}(a,b) = \frac{f(a) + f(b)}{2} - \frac{\alpha(a-b)}{2} \]

where $\alpha = \max(v + c_s)$ and $c_s = \sqrt{\gamma \frac{p}{\rho}}$.

Assuming a simple steady state class of solutions, consider the non-moving equilibria, where $v \equiv 0$, \eqref{eq:1ddensity} has the following (component wise) update:

\[H_\rho^i = \frac{\alpha}{2}\jump{\rho(x_{k+1/2})}\psi_i(1) - \frac{\alpha}{2}\jump{\rho(x_{k-1/2})}\psi_i(-1)\]

where $\jump{f(x)}=f(x^+)-f(x^-)$ denotes the jump in $f$ between left and right states at $x$, and $\average{f(x)} = \frac{f(x^+)+f(x^-)}{2}$ the average of the jump of $f$ at $x$. This shows that the error comes from the jumps in the variable $\rho$ at the interfaces of the control volume $K$.

Similarly, the update function $H_{\rho v}$ for \eqref{eq:1dmomentum} is:

\[H_{\rho v}^i = -\average{ p(x_{k+1/2})}\psi_i(1) + \average{p(x_{k-1/2})} \psi_i(-1) - \sum_j \rho \partial_x\Phi \psi_i(x_j)w_j + \sum_j p \partial_x \psi(x_j)w_j \]

And for \eqref{eq:1denergy}:

\[ H_E^i = \frac{\alpha}{2}\jump{E(x_{k+1/2})}\psi_i(1) - \frac{\alpha}{2}\jump{E(x_{k-1/2})}\psi_i(-1)\]

For the density and energy evolution, the error comes from the jump on the respective variable at the interfaces. For the momentum equation, the error will arise from the split treatment when discretising $\nabla\cdot f(w)$ and $s(w)$, which should exactly cancel out if $w$ is a steady state solution.

Now we consider a general class of steady state solutions, where $v \not\equiv 0$. From \eqref{eq:1ddensity}, follows that $\rho v = \mbox{const}$.

Then, for \eqref{eq:1ddensity} one can write the following update function:

\[H_\rho^i = \sum_j \rho v \partial_x\psi_i(x_j)w_j - \frac{\alpha}{2}\jump{\rho(x_{k+1/2})}\psi_i(1) + \frac{\alpha}{2}\jump{\rho(x_{k-1/2})}\psi_i(-1) - \average{\rho v} \psi_i (1) + \average{\rho v} \psi_i (-1). \]

Rewriting $H_\rho$, using the fact that $\rho v = \mbox{const}$ and that Legendre polynomials have the property $\psi_n(-x)=(-1)^n \psi_n(x)$, one can arrive at: 
\begin{align*}
H_\rho^i &= \frac{\alpha}{2}\jump{\rho(x_{k+1/2})}\psi_i(1) - \frac{\alpha}{2}\jump{\rho(x_{k-1/2})}\psi_i(-1)
\end{align*}
Independently of the order of the polynomial $\psi_n$, the volume integral part cancels out either due to the numerical flux contribution or due to the fact $\rho v$ is constant. 

For \eqref{eq:1dmomentum}:
\begin{align*}
H_{\rho v}^i &= -\average{ \rho v^2 + p (x_{k+1/2})}\psi_i(1)  + \average{\rho v^2 + p (x_{k-1/2})}\psi_i(-1)\\
& - \sum_j \rho \partial_x\phi \psi_i(x_j)w_j + \sum_j (\rho v^2 + p) \partial_x \psi(x_j)w_j 
\end{align*}

And for  \eqref{eq:1denergy}:
\begin{align*}
H_{E}^i &= \frac{\alpha}{2}\jump{E(x_{k+1/2})}\psi_i(1) - \frac{\alpha}{2}\jump{E(x_{k-1/2})}\psi_i(-1) \\ 
&+ \average{ v(E+p)(x_{k+1/2}) }\psi_i(1) - \average{ v(E+p)(x_{k-1/2})}\psi_i(-1) \\
&- \rho v \sum_j  \partial_x\phi \psi_i(x_j)w_j + \sum_j (v(E+p)) \partial_x \psi(x_j)w_j 
\end{align*}
While update $H_{\rho}$ remains unchanged, $H_{\rho v}$ and $H_{E}$ have additional terms from the velocity contribution, and thus one can observe that the error arises from splitting the flux term in surface and volume terms, and the separated treatment of the source term and $\nabla \cdot f(w)$.

\section{Well-balanced RKDG method}
\label{Wellbal}

We now present our implementation of a well balanced method for RKDG. Using the formulation presented in \eqref{eq:dg}, we follow an approach similar to \cite{dedner2001}, 
where we represent the solution of \eqref{eq:main} as a sum of a steady state (or equilibrium) solution $w_{eq}(\bm{x})$ and a perturbation $\delta w(\bm{x},t)$:
\begin{equation*}
    w(\bm{x},t) = w_{eq}(\bm{x}) + \delta w(\bm{x},t) \quad a.e.
\end{equation*}
We note that if \eqref{eq:main} admits a steady state solution $w_{eq}$, the flux-source balance relation holds:
\begin{equation}
\label{eq:fluxbal}
    \nabla \cdot f(w_{eq}(\bm{x})) = s(w_{eq}(\bm{x})).
\end{equation}
And weakly, for a suitable test function $v(\bm{x})$:
\begin{equation}
\label{eq:fluxbalweak}
    \int \nabla \cdot f(w_{eq}(\bm{x})) v(\bm{x}) d\bm{x} = \int s(w_{eq}(\bm{x})) v(\bm{x}) d\bm{x}.
\end{equation}

Subtracting \eqref{eq:fluxbalweak} from \eqref{eq:weakform}, and noting that a state state solution satisfies  $\frac{\partial}{\partial t} w_{eq} = 0$, we can write:

\begin{equation*}
\begin{split}
\frac{\mbox{d}}{\mbox{dt}}\int_K (\delta w(\bm{x},t))v(\bm{x})\mbox{d\bf{x}} =& -\sum_{e\in\partial K}\int_e \delta f(w(\bm{x},t))\cdot n_{e,K} v(\bm{x}) \mbox{d}\Gamma\\ &
\quad + \int_K \delta f(w(\bm{x},t))\cdot \nabla v(\bm{x})\mbox{d\bf{x}}\\ &\quad + \int_K \delta s(w(\bm{x},t)) v(\bm{x})\mbox{d\bf{x}}.
\end{split}\end{equation*}
where we use the following notation:
\begin{enumerate}
    \item{
    $\begin{aligned}[t]
    \int_e \delta f(w(\bm{x},t))\cdot n_{e,K} v(\bm{x}) \mbox{d}\Gamma = \int_e \big( f(w(\bm{x},t)) - f(w_{eq}(\bm{x})) \big) \cdot n_{e,K} v(\bm{x}) \mbox{d}\Gamma
    \end{aligned}$
    }
    \item{
    $\begin{aligned}[t]
    \int_K \delta f(w(\bm{x},t))\cdot \nabla v(\bm{x})\mbox{d\bf{x}} = \int_K (f(w(\bm{x},t))-f(w_{eq}(\bm{x})))\cdot \nabla v(\bm{x})\mbox{d\bf{x}}
    \end{aligned}$
    }
    \item{
    $\begin{aligned}[t]
    \int_K \delta S(w(\bm{x},t)) v(\bm{x})\mbox{d\bf{x}} = \int_K (s(w(\bm{x},t))-s(w_{eq}(\bm{x}))) v(\bm{x})\mbox{d\bf{x}}
    \end{aligned}$
    }
\end{enumerate}

Note again that there is an ambiguity in the definition of the flux $f(w(\bm{x},t))\cdot n_{e,K}$  since $w$ can be multi-valued at the cell interface. To overcome this inconsistency, the ambiguous term is here again replaced by a single-valued numerical flux $h_{e,K}(\bm{x},t)$ computed using the Lax Friedrich Riemann solver.

 
Let our numerical solution be represented as:

\[ w_{num}(\bm{x},t) = w_{eq}(\bm{x}) + \delta w_h (\bm{x}, t), \]

where  $\delta w_h \in V_h(K)$. Furthermore, we approximate the integrals with a quadrature, which yields the following well balanced DG numerical scheme:
\begin{equation*}
\begin{split}
 \delta w_h(t=0)&=P_{V_h}(\delta w_0) \\
 \frac{\mbox{d}}{\mbox{dt}}\int_K \delta w_h(\bm{x},t)v_h(\bm{x})\mbox{d\bf{x}} &= -\sum_{e\in\partial K}\sum_{i=0}^L \delta f_{e,K}(w_{num}(\bm{x}_i,t))v_h(\bm{x}_i) \omega_i |e| \\
 &\quad + \sum_{j=0}^M \delta f(w_{num}(\bm{x}_j,t))\cdot \nabla v_h(\bm{x}_j) \omega_j |K| \\
 &\quad+ \sum_{j=0}^M \delta s(w_{num}(\bm{x}_j,t)) v_h(\bm{x}_j) \omega_j |K| \qquad \forall v_h(\bm{x})\in V(K), \forall K \in T_h.
 \end{split}
\end{equation*}

Note that this reformulation is only suitable for problems where the solution $w$ is close \emph{enough} to the prescribed steady state solution $w_{eq}$. In fact, if the initial condition is exactly equal to the steady state solution ($w_0 = w_{eq}$), the scheme will capture the equilibrium solution exactly. If the initial condition is close to the steady state solution, this scheme is able to evolve the perturbation without being dominated by the truncation error on the steady state solution. However, if the initial condition is very far from the adopted steady state, this scheme might not be suitable and the traditional RKDG scheme will be more robust. If this is the case, setting the steady state $w_{eq}$ to 0 and the perturbation $\delta w$ to the full solution, one simply recovers the traditional RKDG scheme \cite{cs5}.
\section{Numerical Experiments}
\label{sec:results}
In this section, several benchmark problems will be introduced. These will be the basis of our discussion in Section \ref{sec:conclusion}. An introduction and description of code used to perform the numerical experiments can be found in \cite{david2018}. Additional results can be found in appendix~\ref{ap:extraresults}.

\subsection{Error estimate}

The empirical error estimates are calculated using the $\mathcal{L}_1$-error norm:
\begin{equation*}
    ||w_h(\bm{x}) - w(\bm{x}) ||_{1} = \int_D |w_h(\bm{x}) - w(\bm{x}) | d\bm{x}. 
\end{equation*}
It is shown in \cite{convergencel1} that a convergence rate of $N_p+1$ for a $N_p$ degree polynomial approximation of the solution in $\mathcal{L}_1$-error norm is expected for smooth enough functions. This quantity is computed with a numerical quadrature and computed the following manner:

\begin{equation}
    ||w_h(\bm{x}) - w(\bm{x}) ||_{1} \approx  \sum_{K \in D} \sum_{i = 0}^{M} \sum_{j = 0}^{M}  \mid w_h(\nu^K(x_i,y_j)) - w(\nu^K(x_i,y_j)) \mid \omega _i \omega_j \frac{\Delta x \Delta y}{4},
 \end{equation}
where $\{x_i,y_j\}_{i,j=0}^M$ are Gauss Legendre quadrature points, $\{\omega_i, \omega_j\}_{i,j=0}^M$ the corresponding weights and $\nu^K(\chi,\upsilon)$ a linear transformation mapping element K to the canonical element $[-1,1]\times[-1,1]$,

\[ \nu^K(\chi,\upsilon) = (x_l - \chi\frac{\Delta x}{2}, y_l - \upsilon\frac{\Delta y}{2} ), \]

and $(x_l,y_l)$ the center of element $K$.

\subsection{Well-balanced property}
\label{subsec1}
In this section, the well-balanced property of the schemes is evaluated. To this end, we first evolve a hydrostatic equilibrium solution. What should be observed, in this case, is that the solution does not change for any time $T>0$. However, due to the failure of perfectly balancing the discrete version of $\nabla\cdot f(w)$ and $s(w)$, the state at some time $T$ might deviate from the initial condition. We then solve for the propagation of perturbations of the equilibrium solution, that we call here {\it waves}, using various amplitudes, and measure whether the schemes can capture these perturbations without being affected by the truncation errors of the equilibrium solution. In the last set of test cases, we evaluate the quality of our schemes using a dynamical equilibrium state, meaning that the velocity $\mathbf{v}=(v_x,v_y)$ is non-zero for the steady state solution.

\subsubsection{Hydrostatic equilibrium}
\paragraph{{1-dimensional case}}\;
Considering an ideal gas $\gamma = 1.4$ and a linear gravitational potential $\Phi_x = gx$, we are interested in preserving the following isothermal equilibrium state:
\begin{equation}
\begin{split}
\label{eq:1dhydrostatic}
    \rho_{eq}(x) &= \rho_0 \exp \bigg(-\frac{\rho_0 g}{p_0}x \bigg) \\
    u_{eq}(x) &= 0\\
    p_{eq}(x) &= p_0 \exp \bigg(-\frac{\rho_0 g}{p_0}x\bigg) \\
\end{split}
\end{equation}
with $\rho_0 = 1.0$, $p_0 = 1.0$ and $g=1.0$.

Because we are interested in preserving the equilibrium state, we impose the boundary condition as the extension of the domain in $\partial D$, as follows:
\begin{equation}
\label{bcs}
\rho(\bm{x})\rvert_{\bm{x}\in \partial D} = \rho_{eq}(\bm{x}) \quad v_x(\bm{x})\rvert_{\bm{x}\in \partial D} = v_{x,eq}(\bm{x}) \quad  v_y(\bm{x})\rvert_{\bm{x}\in \partial D} = v_{x,eq}(\bm{x}) \quad p(\bm{x})\rvert_{\bm{x}\in \partial D}  = p_{eq}(\bm{x}),
\end{equation}
where $\bm{x} = (x) $ in 1-dimension and $\bm{x} = (x,y)$ in 2-dimensions.

\begin{figure}
    \centering
     \includegraphics[width=0.9\textwidth]{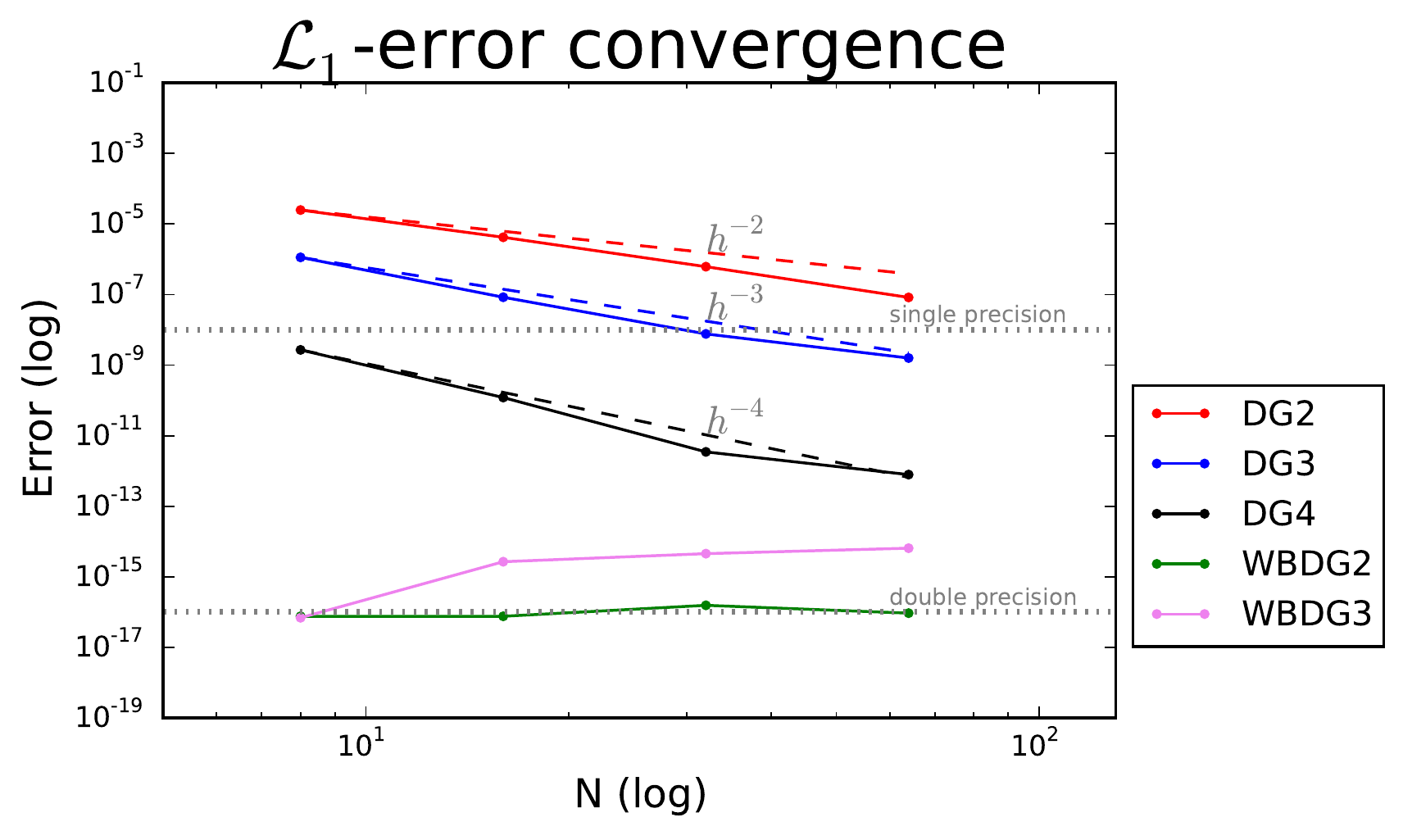}
     \caption{$\mathcal{L}_1$ error convergence for the 1-dimensional hydrostatic test case \eqref{eq:1dhydrostatic}.}
    \label{fig:hydroeq_1d}
\end{figure}

The numerical errors for the density are shown in figure~\ref{fig:hydroeq_1d} for the following resolutions $N = 8,~16,~32,~64$ at time $T = 10.0$ for the second and third order well-balanced scheme (WBDG2 and WBDG3, respectively) and the traditional discontinuous Galerkin method with orders 2, 3 and 4 (DG2, DG3, DG4 respectively). One can observe that by increasing the resolution or the order, the truncation error can be reduced, even for long time evolution. Furthermore, at $N = 64$ at $4^{th}$ order, we reach a similar absolute error as in our well-balanced methods. 

\begin{figure}
    \centering
     \includegraphics[width=0.9\textwidth]{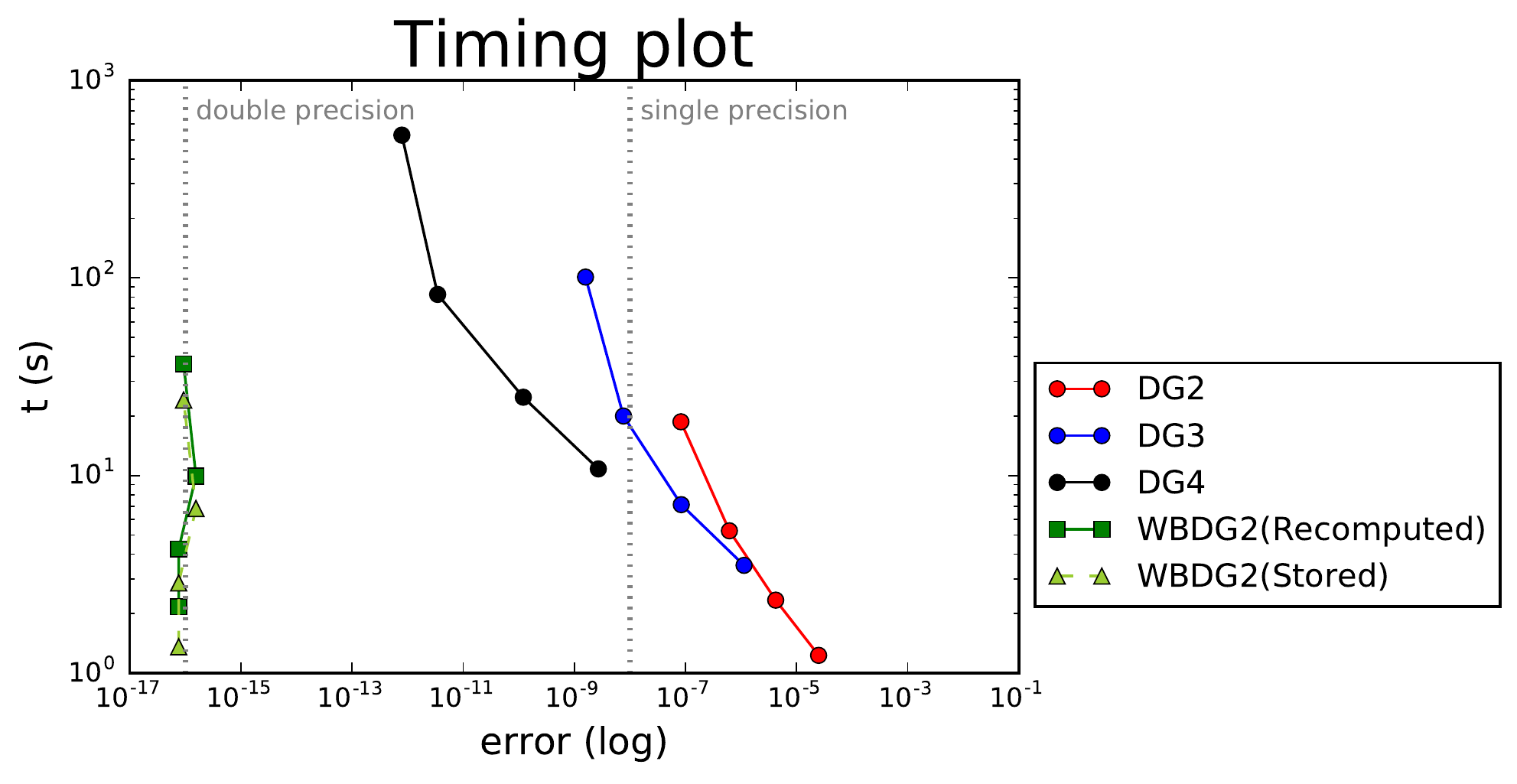}
     \caption{Total time to achieve a particular accuracy for the 1-dimensional hydrostatic test case \eqref{eq:1dhydrostatic}}
         \label{fig:timing_nomem}
\end{figure}

It is important to stress that the well-balanced scheme requires either the storage of additional arrays or requires to perform many additional computation every time step. Indeed, at each Runge-Kutta timestep, we can either recompute face nodal values or or we can store the equilibrium solution once and for all, requiring $\mathcal{O}((4+m)N_x N_y m)$ of additional memory, where $N_x$ and $N_y$ denotes the number of cells in $x-,y-$ direction, respectively, and $m$ the order of the method\footnote{Further optimization is possible, by storing the resulting volume/surface integral for each cell, further reducing the necessary storage to $\mathcal{O}(N_x N_y)$}. Shown in figure~\ref{fig:timing_nomem}, we show the total time it takes to run the 1-dimensional hydrostatic equilibrium test case \eqref{eq:1dhydrostatic} when performing the well-balanced reconstruction (denoted as WBDG2(Rec)) versus precomputing and storing the equilibrium variables (denoted as WBDG2(Mem)), compared to the traditional discontinuous Galerkin methods with order 2, 3 and 4 (denoted as DG2, DG3, DG4 respectively).

A perturbation is now added to the pressure state of the equilibrium solution described in \eqref{eq:1dhydrostatic}, as shown below:
\begin{equation}
\label{eq:pert1d}
\begin{split}
    p(x,t=0) &= p_{eq}(x) + \eta \exp \bigg(-\frac{\rho_0 g}{p_0}\frac{(x-0.5)^2}{0.01}\bigg) \\
\end{split}
\end{equation}

\begin{figure}
    \centering
    \begin{subfigure}[b]{0.49\textwidth}
        \includegraphics[width=\textwidth]{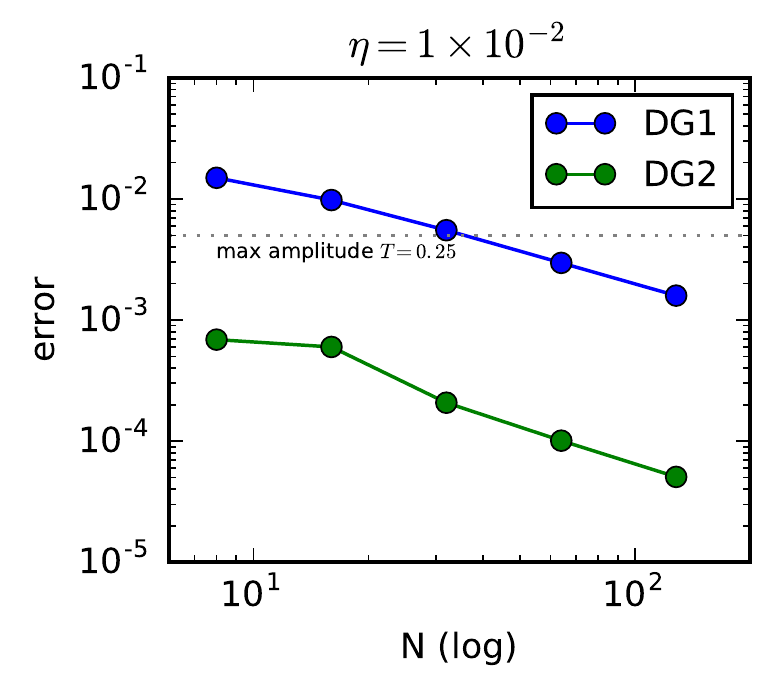}
            \end{subfigure}
                \begin{subfigure}[b]{0.49\textwidth}
        \includegraphics[width=\textwidth]{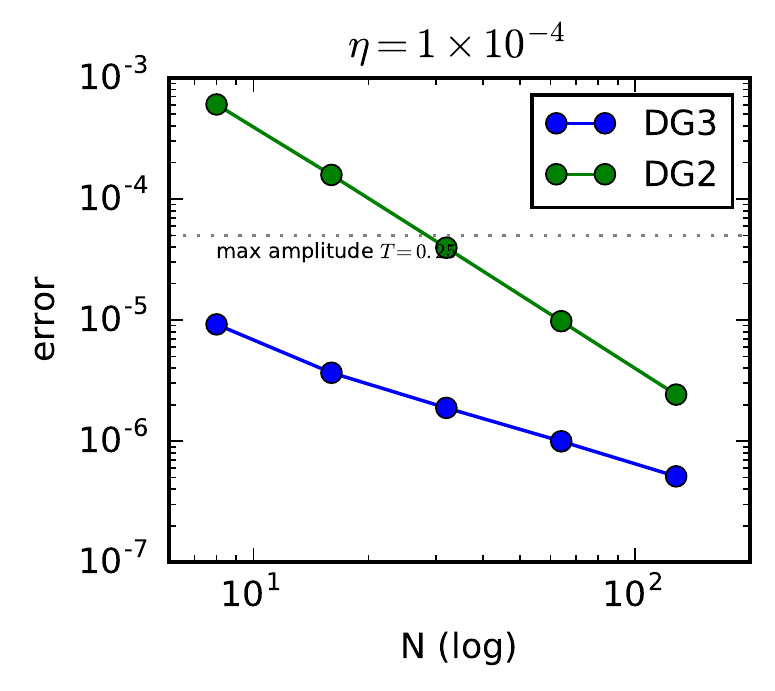}
    \end{subfigure}\\
    \begin{subfigure}[b]{0.49\textwidth}
        \includegraphics[width=\textwidth]{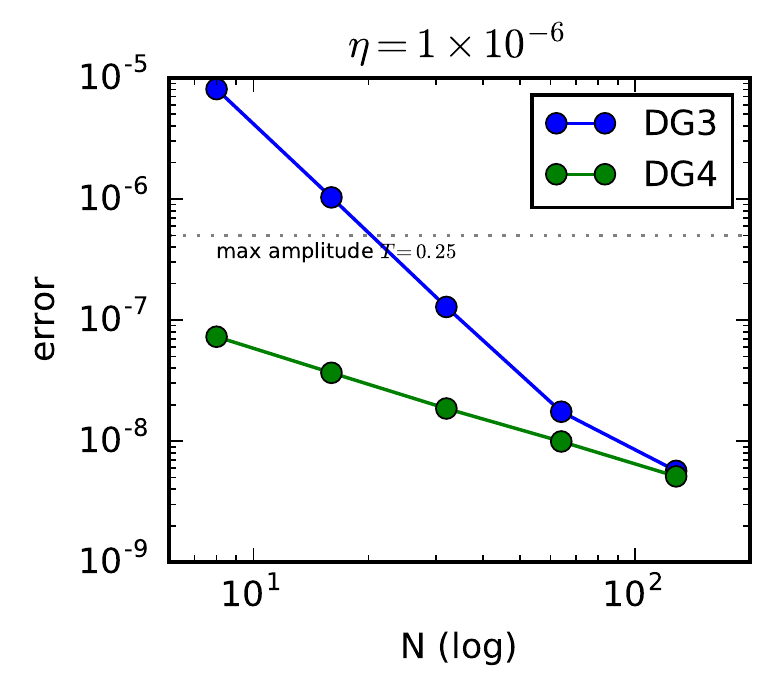}
    \end{subfigure}
    \begin{subfigure}[b]{0.49\textwidth}
        \includegraphics[width=\textwidth]{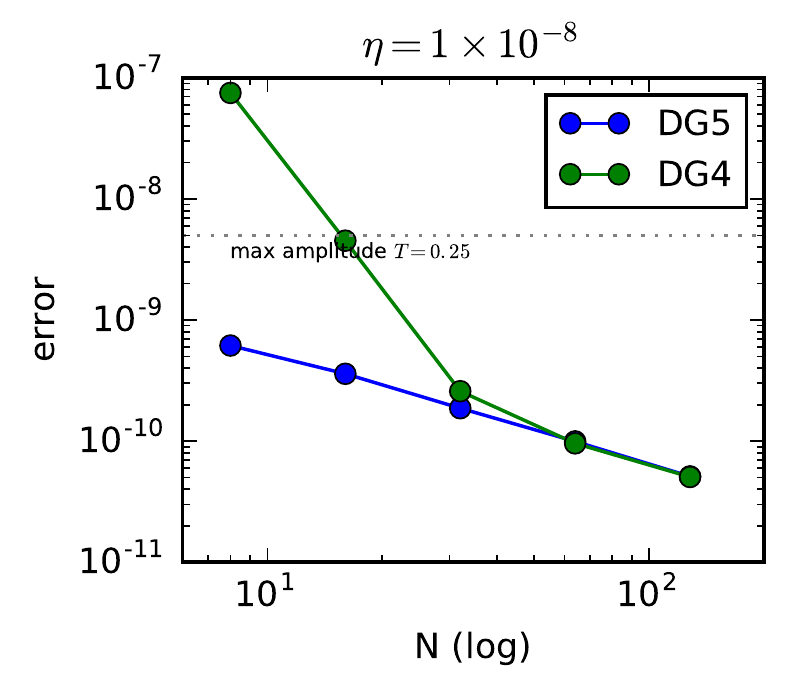}
    \end{subfigure}
    \caption{Convergence of the non well-balanced methods for initial conditions \eqref{eq:pert1d} with perturbation sizes of $\eta = 1\times10^{-2},~1\times10^{-4},~1\times10^{-6}$ and $1\times10^{-8}$ respectively.}
    \label{fig:conv1dpulse}
\end{figure}

\begin{figure}
    \centering
    \begin{subfigure}[b]{0.49\textwidth}
        \includegraphics[width=\textwidth]{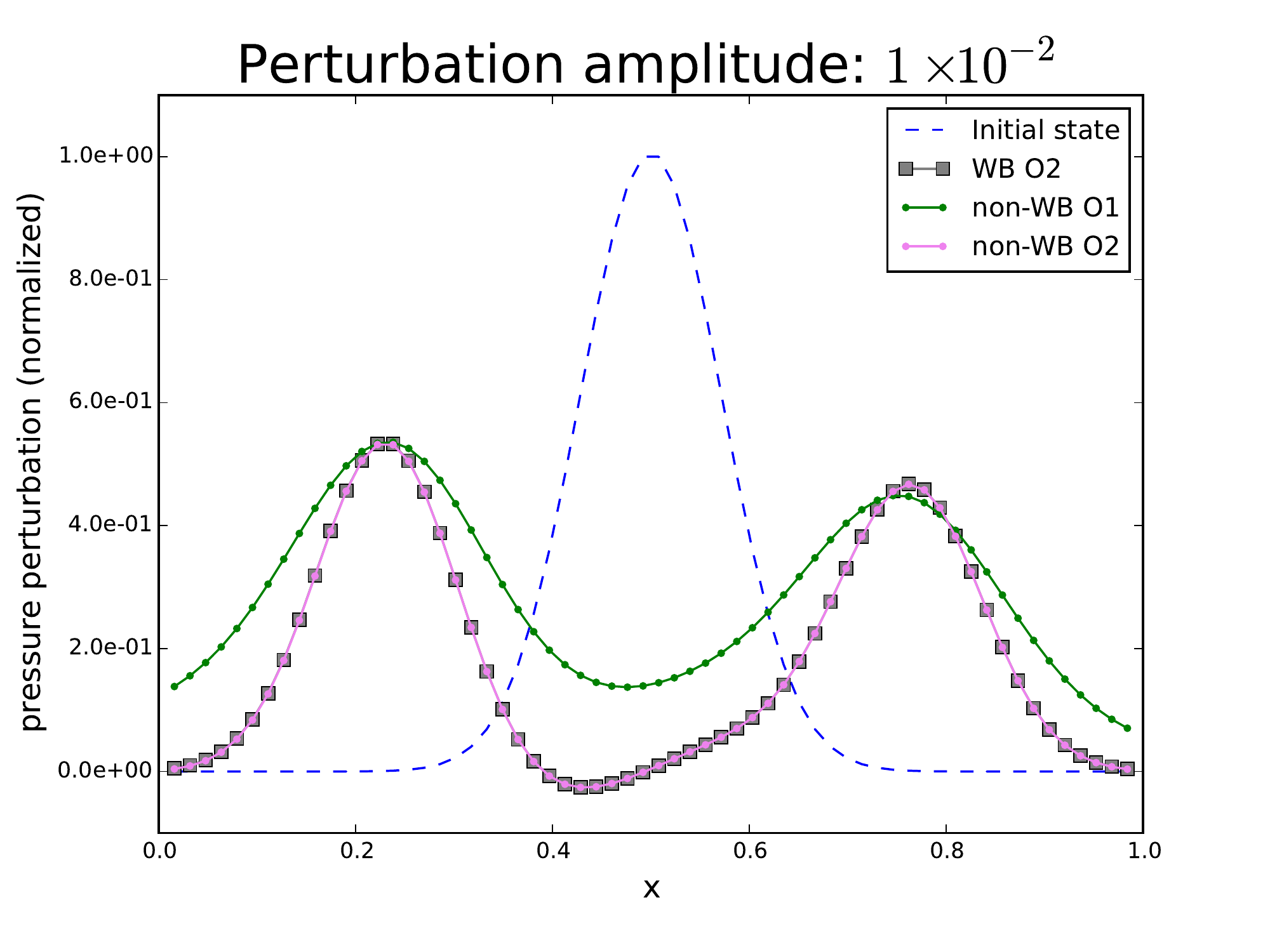}
            \end{subfigure}
                \begin{subfigure}[b]{0.49\textwidth}
        \includegraphics[width=\textwidth]{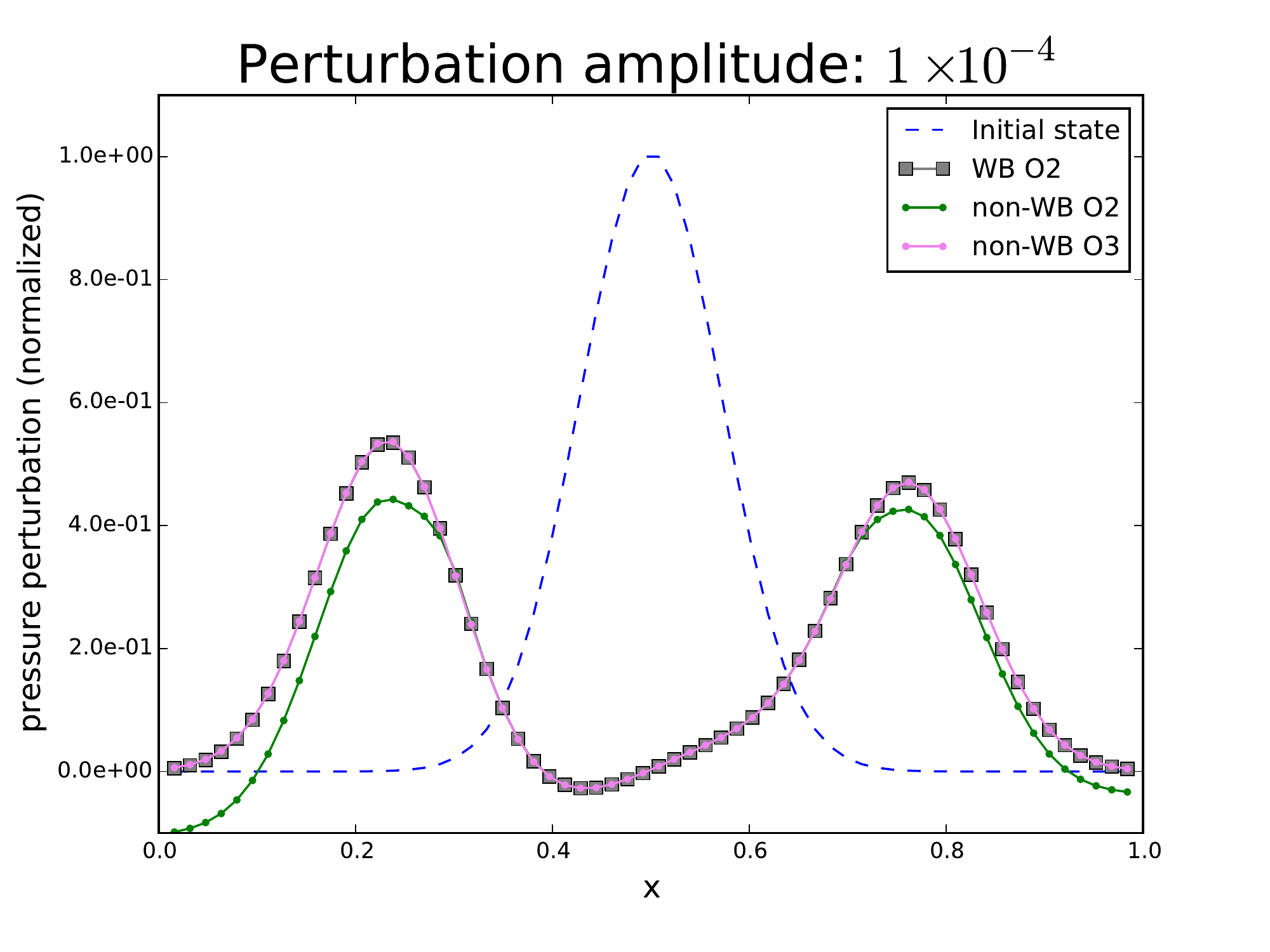}
    \end{subfigure}\\
    \begin{subfigure}[b]{0.49\textwidth}
        \includegraphics[width=\textwidth]{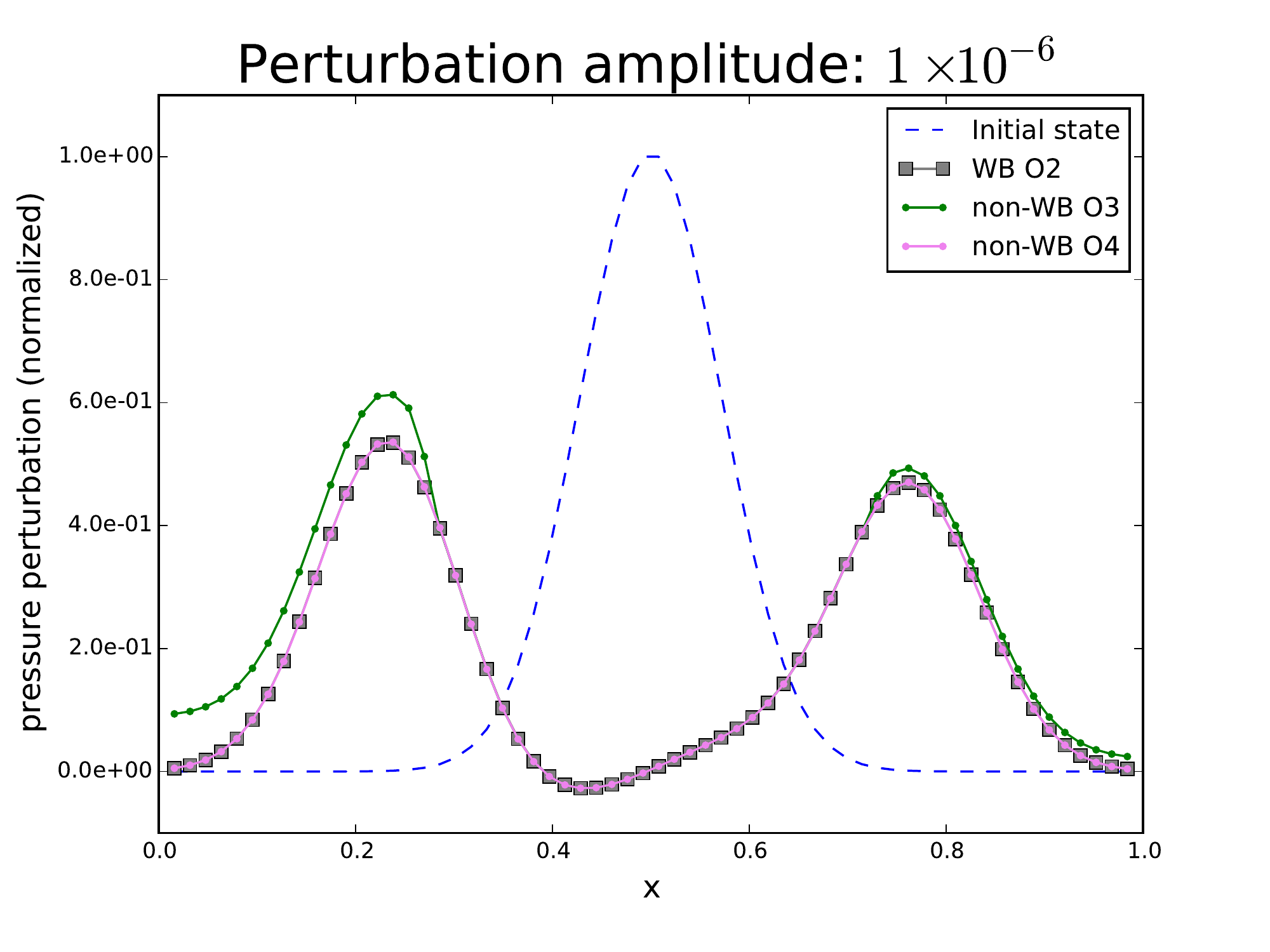}
    \end{subfigure}
    \begin{subfigure}[b]{0.49\textwidth}
        \includegraphics[width=\textwidth]{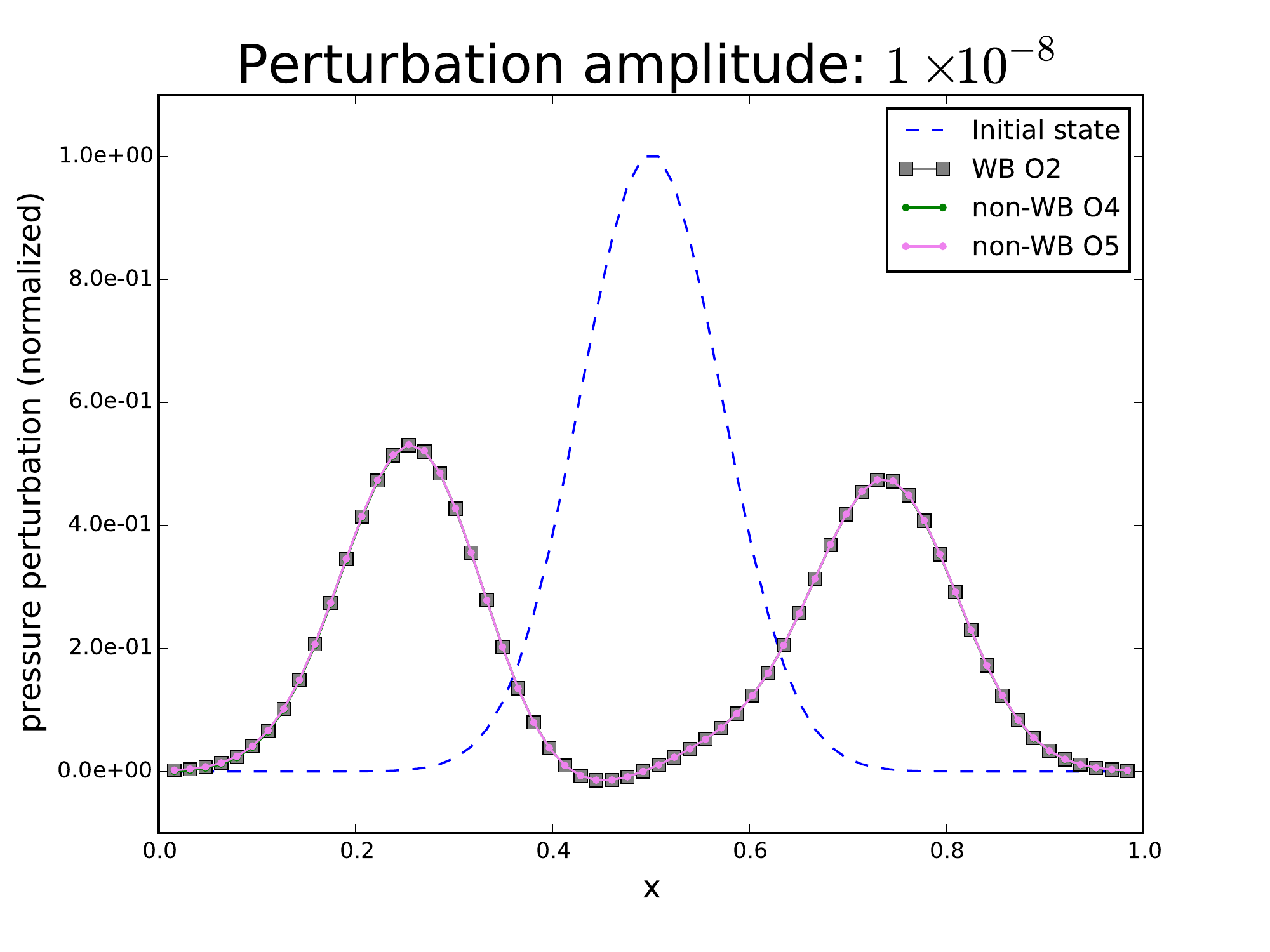}
    \end{subfigure}
    \caption{Non well-balanced method versus well-balanced method for hydrostatic equilibrium with varying amplitude perturbation on the pressure field for initial conditions \eqref{eq:pert1d}.}
    \label{fig:1dpulses}
\end{figure}

The initial condition \eqref{eq:pert1d} is run until $T=0.25$ with different pulse amplitudes: $\eta = 1\times10^{-2},~1\times10^{-4},~1\times10^{-6}$ and $1\times10^{-8}$. In figure~\ref{fig:conv1dpulse} we show the pointwise $L_1$ error of between the solution attained with different orders of the non-well balanced discontinuous Galerkin scheme and a high resolution solution which captures the pulse. We note that if the error is larger than the perturbation size, then it is clear that a particular combination of order and resolution is not enough to capture the perturbation. A qualitative depiction of this is shown in figure~\ref{fig:1dpulses}, for a fixed grid-size of $N = 64$ and different orders. Furthermore, we note that for $\eta=1\times10^{-2}$, the difference between a second-order well balanced and a second-order non-well-balanced scheme is impossible to see. However, when the perturbation's amplitude $\eta$ decreases below the truncation error of the scheme, the wave is no longer well captured. As shown in figure~\ref{fig:1dpulses}, the error for a non-well-balanced scheme can be reduced by increasing the order of the scheme or the resolution of the grid, effectively reducing the approximation error. Note that for the well-balanced method, we always capture the correct wave solution. The time to solution for the experiment with perturbation size $1\times 10^{-8}$ is shown in table~\ref{tab:hydrostatic1-8}. Additional times to solution can be found in appendix~\ref{ap:extraresults}.

\begin{table}[h]
\centering
\caption{Time to solution for initial conditions \eqref{eq:pert1d} for $\eta = 1\times 10^{-8}$ in seconds (s).}
\label{tab:hydrostatic1-8}
\begin{tabular}{||c | c c c c ||} 
 \hline
 N & DG4& DG5 & WBDG2 & WBDG3 \\ [0.5ex] 
 \hline\hline
 8 & 0.37 & 0.56 & 0.05 & 0.14\\ 
 16 & 0.62 & 1.28 & 0.10 & 0.29\\ 
 32 & 2.05 & 4.79 & 0.20 & 0.69\\ 
 64 & 12.8 & 32.4 & 0.62 & 3.25\\
 128 & 103 & 270 & 3.84 & 24.1\\   [1ex] 
 \hline
\end{tabular}
\end{table}

\paragraph{{2-dimensional case}}
We consider an ideal gas $\gamma = 1.4$ and a linear gravitational potential $\Phi = g(x + y)$. We are interested in preserving the following isothermal equilibrium state on a unit square domain $\bm{x}\in[0,1]\times[0,1]$:
\begin{equation}
\label{eq:2dhydrostatic}
\begin{split}
    \rho_{eq}(x,y) &= \rho_0 \exp \bigg(-\frac{\rho_0 g}{p_0}(x+y) \bigg) \\
    u_{eq}(x,y) &= 0\\
    v_{eq}(x,y) &= 0\\
    p_{eq}(x,y) &= p_0 \exp \bigg(-\frac{\rho_0 g}{p_0}(x+y) \bigg) \\
\end{split}
\end{equation}
with $\rho_0 = 1$, $p_0 = 1$ and $g=1$.

\begin{figure}
    \centering
     \includegraphics[width=0.8\textwidth]{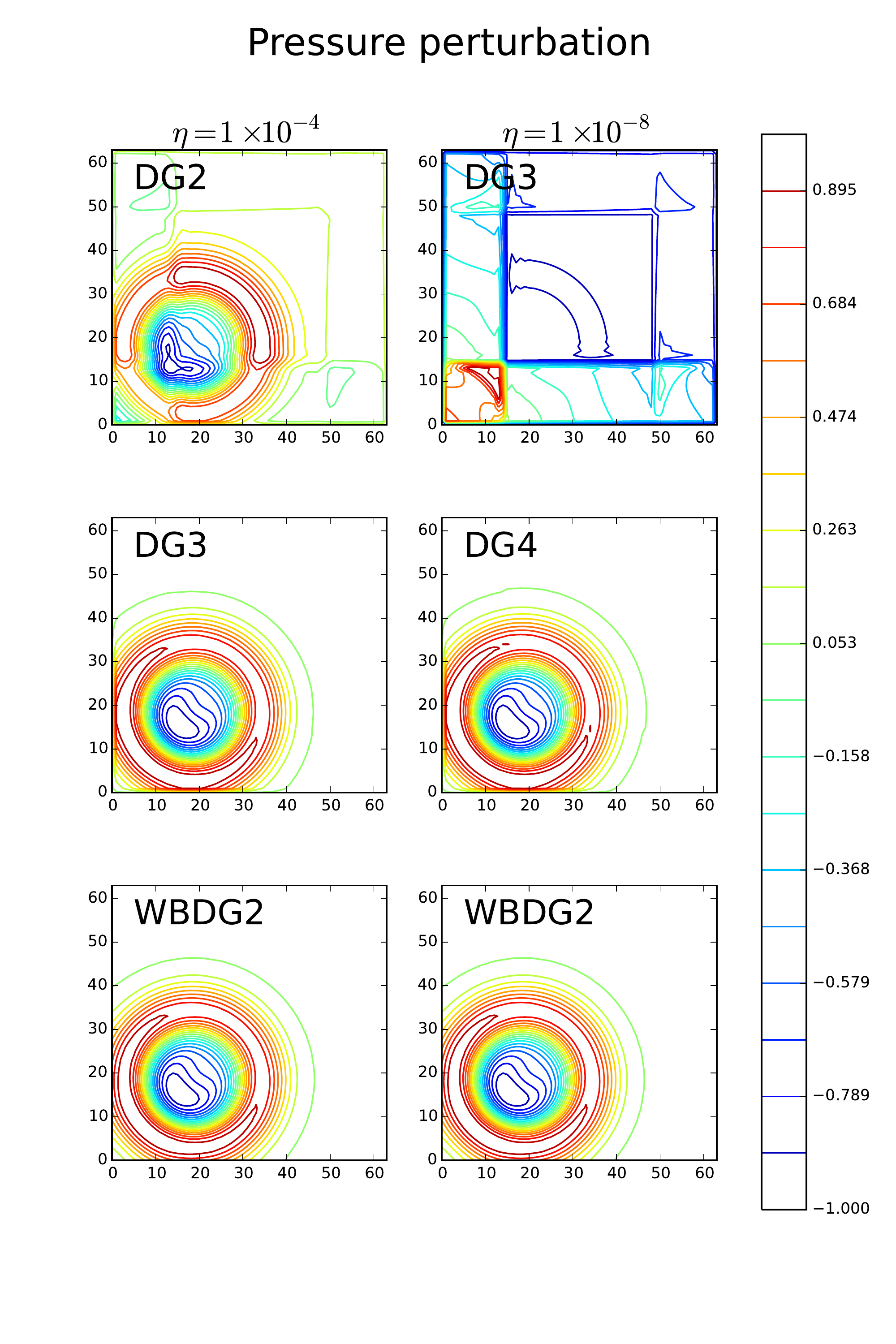}
     \caption{Non well-balanced method vs well-balanced method for hydrostatic equilibrium with varying amplitude perturbation on the pressure field for initial conditions \eqref{eq:pert2d}.}\label{fig:hydro_2d}
\end{figure}

The numerical errors for the pressure are reported in figure~\ref{fig:hydro_2d} for the following resolutions $N = 8,~16,~32,~64$, evaluated at final time $T = 10.0$. Similarly to the 1-dimensional case, one can observe that the truncation error can be reduced again by increasing the number of cells or the order, as expected.

\begin{figure}
    \centering
     \includegraphics[width=0.9\textwidth]{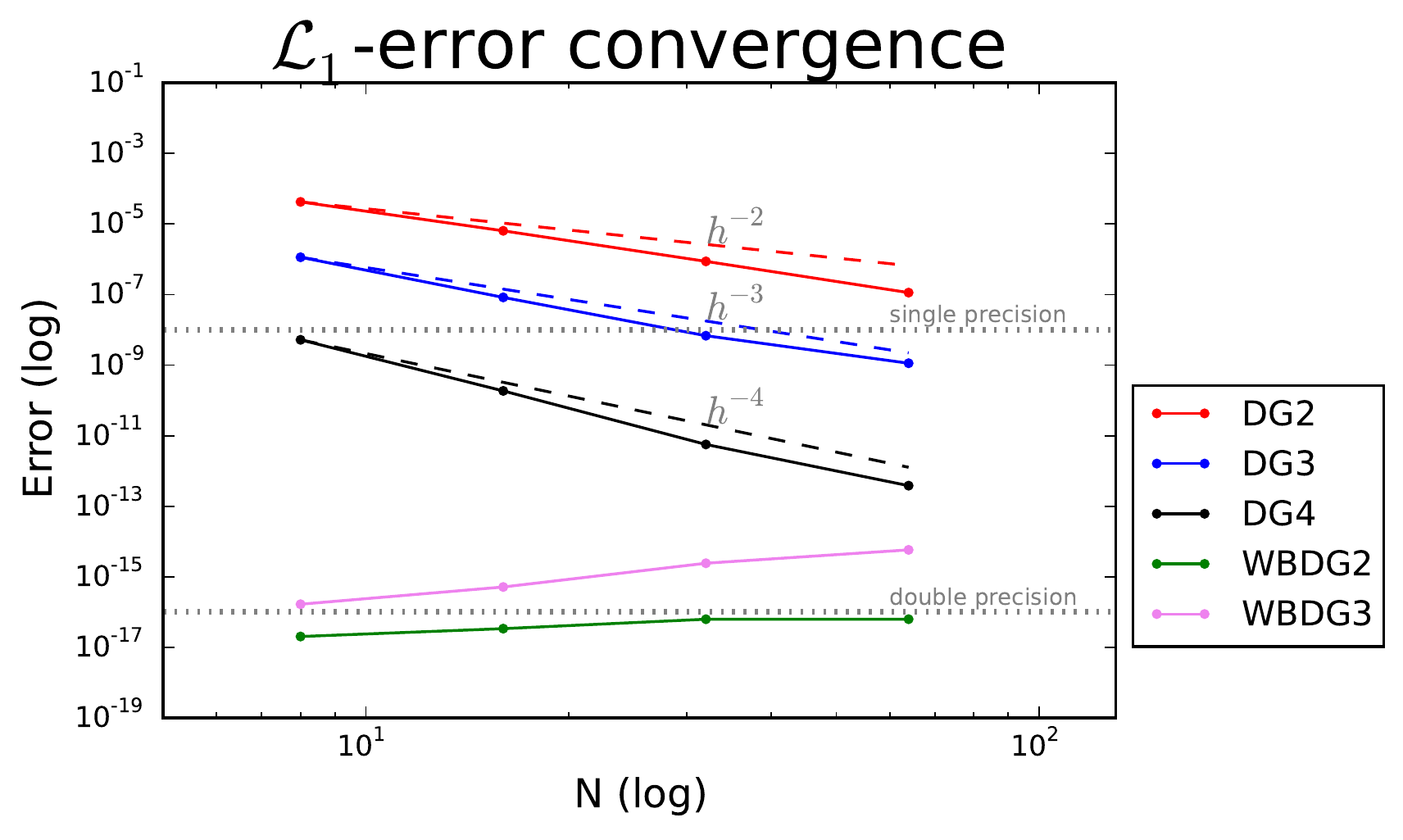}
     \caption{$\mathcal{L}_1$ error convergence for the 2-dimensional hydrostatic test case \eqref{eq:2dhydrostatic}.}
\end{figure}

Again, as in 1-dimension, a perturbation is added to the pressure state of the isothermal equilibrium solution: 
\begin{equation}
\label{eq:pert2d}
\begin{split}
    p(x,y,t=0) &= p_{eq}(x,y) + \eta \exp \bigg(-\frac{\rho_0 g}{p_0}\bigg(\frac{(x-0.3)^2+(y-0.3)^2}{0.01}\bigg) \bigg) \\
\end{split}
\end{equation}

The initial condition \eqref{eq:pert2d} is run with different pulse amplitudes: $\eta = ~1\times10^{-4}$ and $1\times10^{-8}$. The results are shown in figure~\ref{fig:hydro_2d}. Again, we observe that by increasing the order, we can resolve for small perturbations, but we have to choose the resolution carefully to guarantee that the pulse is captured accurately. As before, the well-balanced methods capture the wave solution correctly, even with a second-order scheme. Further analysis, such as the pointwise $L_1$  error of between the solution attained with different orders of the non-well balanced discontinuous Galerkin scheme and a high resolution solution and simulation time to solution can be found in in appendix~\ref{ap:extraresults}.

\subsubsection{Non-hydrostatic steady state}
\begin{figure}
    \centering
         \includegraphics[width=0.9\textwidth]{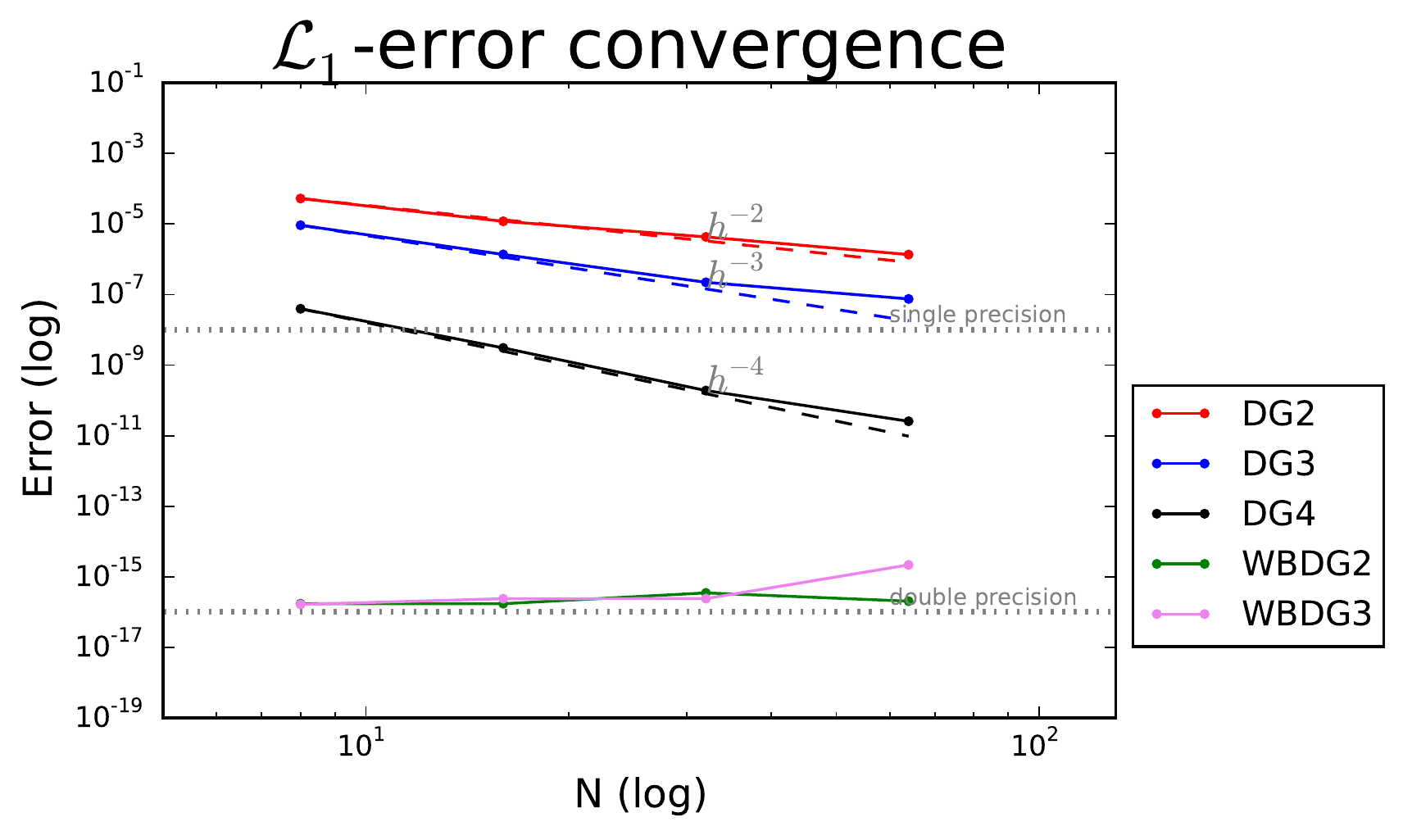}
         \caption{$\mathcal{L}_1$-error convergence for the 1-dimensional dynamic test case \eqref{eq:1ddynamic}.}
         \label{fig:dymeq_1d}
\end{figure}

\paragraph{1-dimensional case}
We consider the manufactured example\footnote{For details of this initial condition, refer to appendix~\ref{Appendix1}.} of an ideal steady gas $\gamma = 1.4$ with a nonzero velocity field and a gravitational field which balances the flux term exactly. We are interested in preserving the following \emph{moving} equilibrium state:
\begin{equation}
\label{eq:1ddynamic}
\begin{split}
    \rho_{eq}(x) &= \rho_0 \exp \bigg(-\frac{\rho_0 g}{p_0}x \bigg) \\
    u_{eq}(x) &= \exp(x)\\
    p_{eq}(x) &= \exp \bigg(-\frac{\rho_0 g}{p_0}x \bigg)^{\gamma} \\
\end{split}
\end{equation}
with $\rho_0 = 1$, $p_0 = 1$ and a non linear potential $\phi = \exp(x) (-\exp(x) + \gamma \exp(-\gamma x)) $. The boundary values are imposed as in \eqref{bcs}. The results are shown in figure~\ref{fig:dymeq_1d} for $T=10.0$.

Now, just as in the hydrostatic equilibrium case \eqref{eq:pert1d}, a perturbation is added to the pressure field:
\begin{equation}
\label{eq:pert1ddym}
\begin{split}
    p(x,t=0) &= p_{eq}(x) + \eta \exp \bigg(-\frac{\rho_0 g}{p_0} \frac{(x-0.3)^2}{0.01} \bigg) \\
\end{split}
\end{equation}

\begin{figure}
    \centering
    \begin{subfigure}[b]{0.49\textwidth}
        \includegraphics[width=\textwidth]{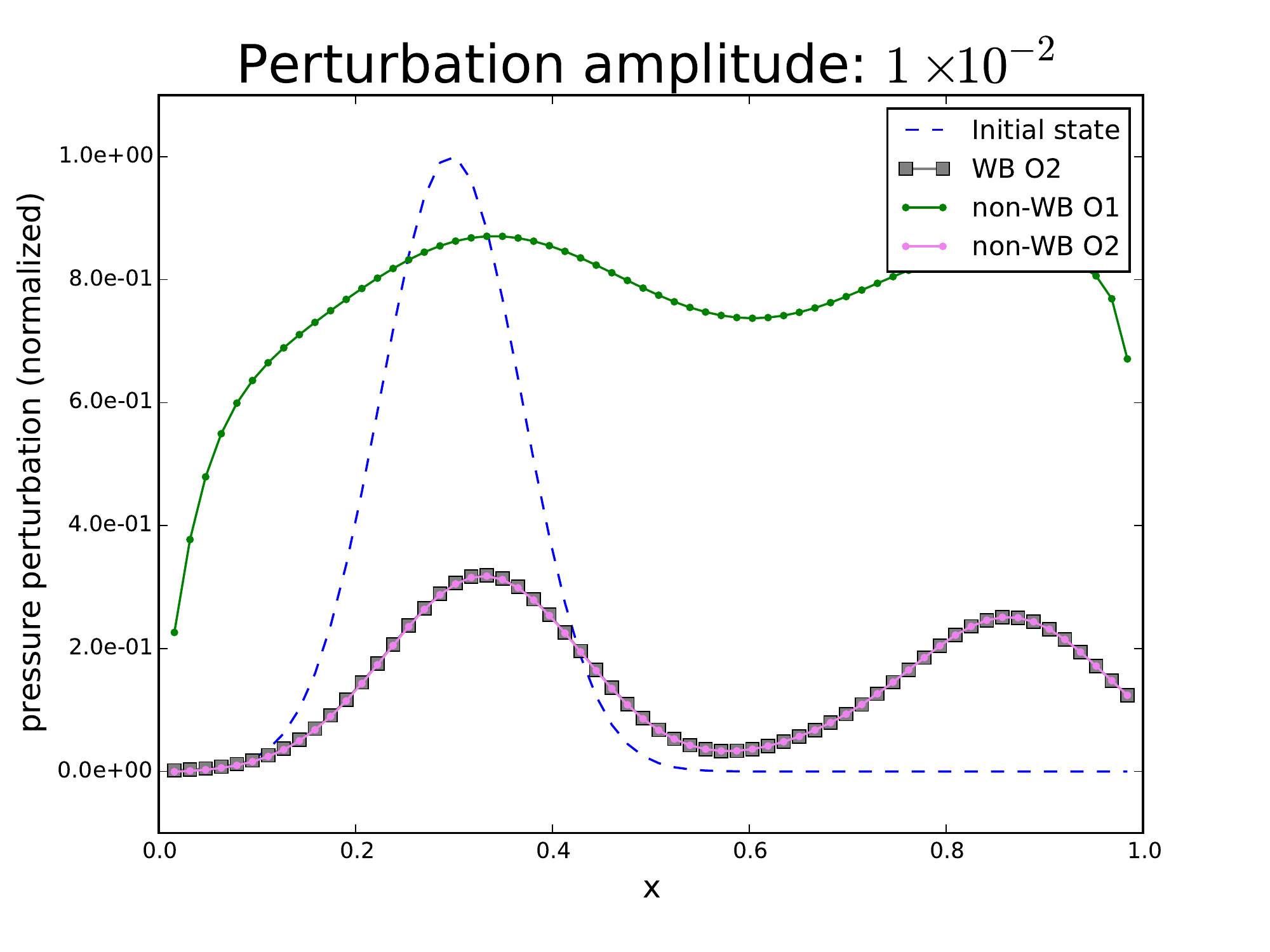}
    \end{subfigure}
    \begin{subfigure}[b]{0.49\textwidth}
        \includegraphics[width=\textwidth]{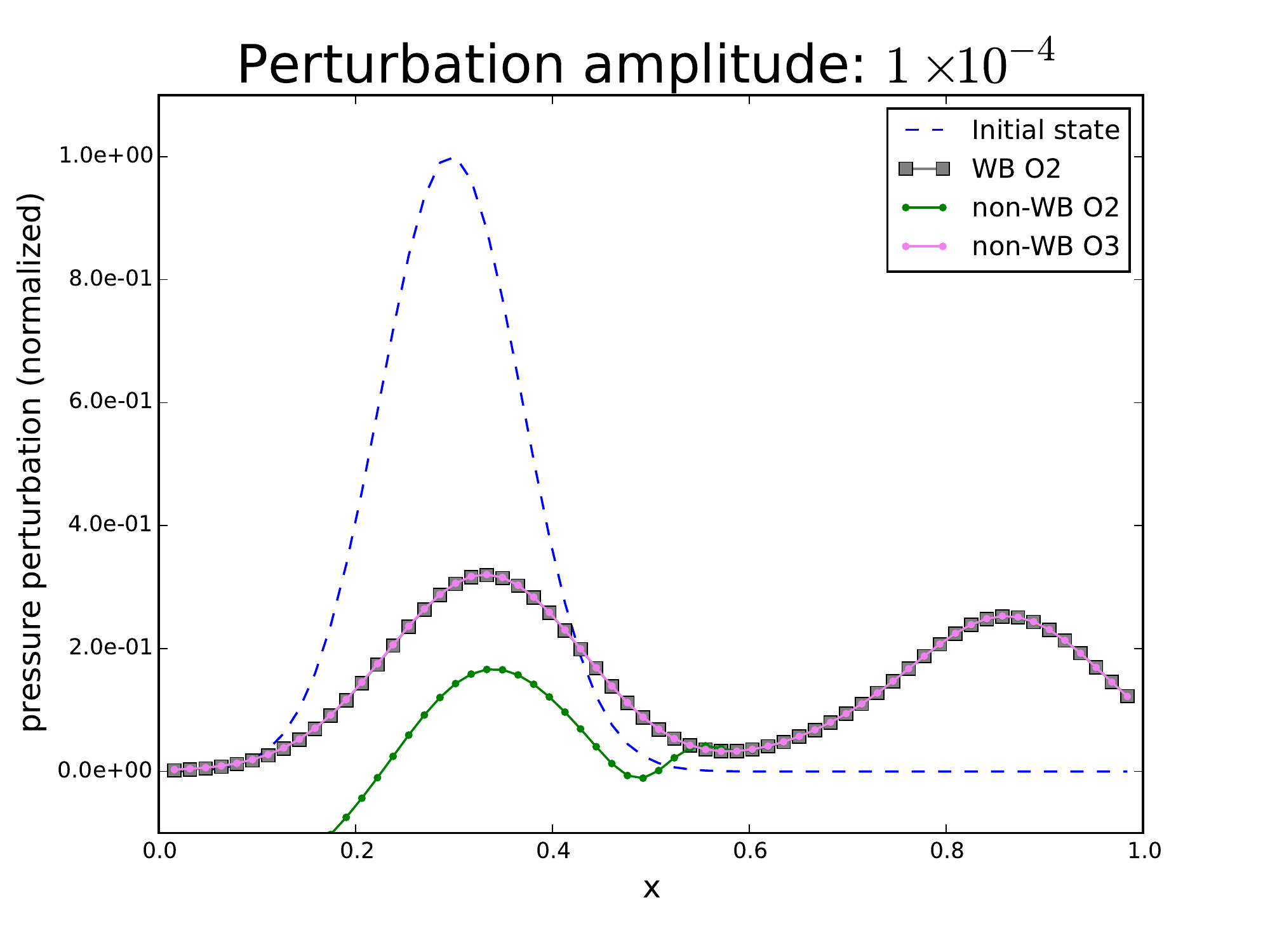}
    \end{subfigure}\\
    \begin{subfigure}[b]{0.49\textwidth}
        \includegraphics[width=\textwidth]{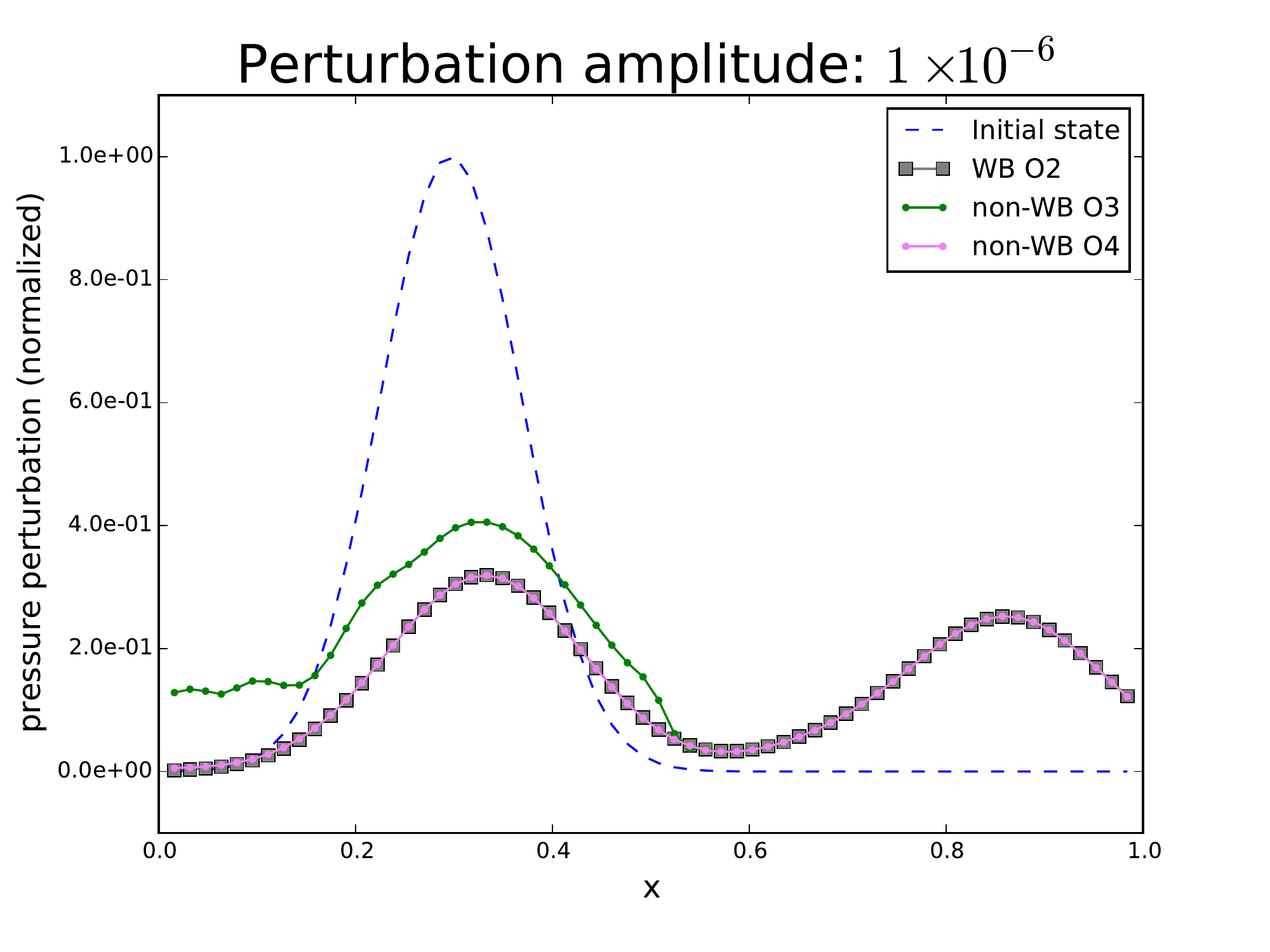}
    \end{subfigure}
    \begin{subfigure}[b]{0.49\textwidth}
     \includegraphics[width=\textwidth]{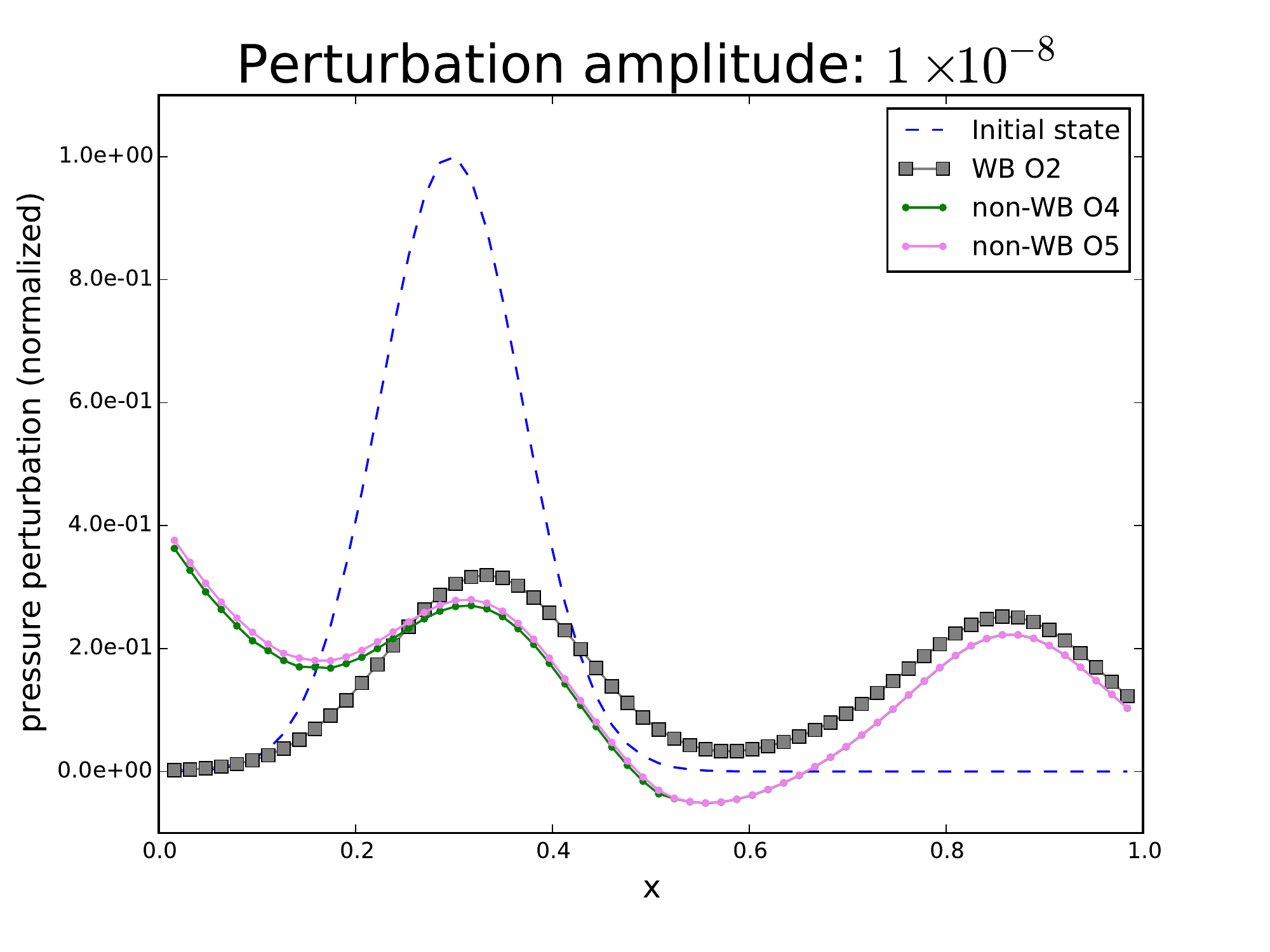}
    \end{subfigure}
    \caption{Non well-balanced method vs well-balanced method for dynamic equilibrium with varying amplitude perturbation on the pressure field as described in \eqref{eq:pert1ddym}.}
    \label{fig:1ddympulses}
\end{figure}

We run the numerical experiment with different pulse amplitudes: $\eta =~1\times10^{-2},~1\times10^{-4},~1\times10^{-6}$ and $1\times10^{-8}$. The results are shown in figure~\ref{fig:1ddympulses}.  Our conclusions remain the same as for the hydrostatic case: for non-well-balanced methods, only a very high order scheme  can capture the low amplitude wave correctly. It appears from figure~\ref{fig:1ddympulses} that the largest truncation error arises from the left boundary and propagates in the direction of the flow. On the contrary, our second-order, well balance method can deal with vanishingly small amplitude waves. Further analysis, such as the pointwise $L_1$  error of between the solution attained with different orders of the non-well balanced discontinuous Galerkin scheme and a high resolution solution and simulation time to solution can be found in in appendix~\ref{ap:extraresults}.

\paragraph{2-dimensional case}
\subparagraph{Modified steady vortex}

\begin{figure}
    \centering
         \includegraphics[width=0.9\textwidth]{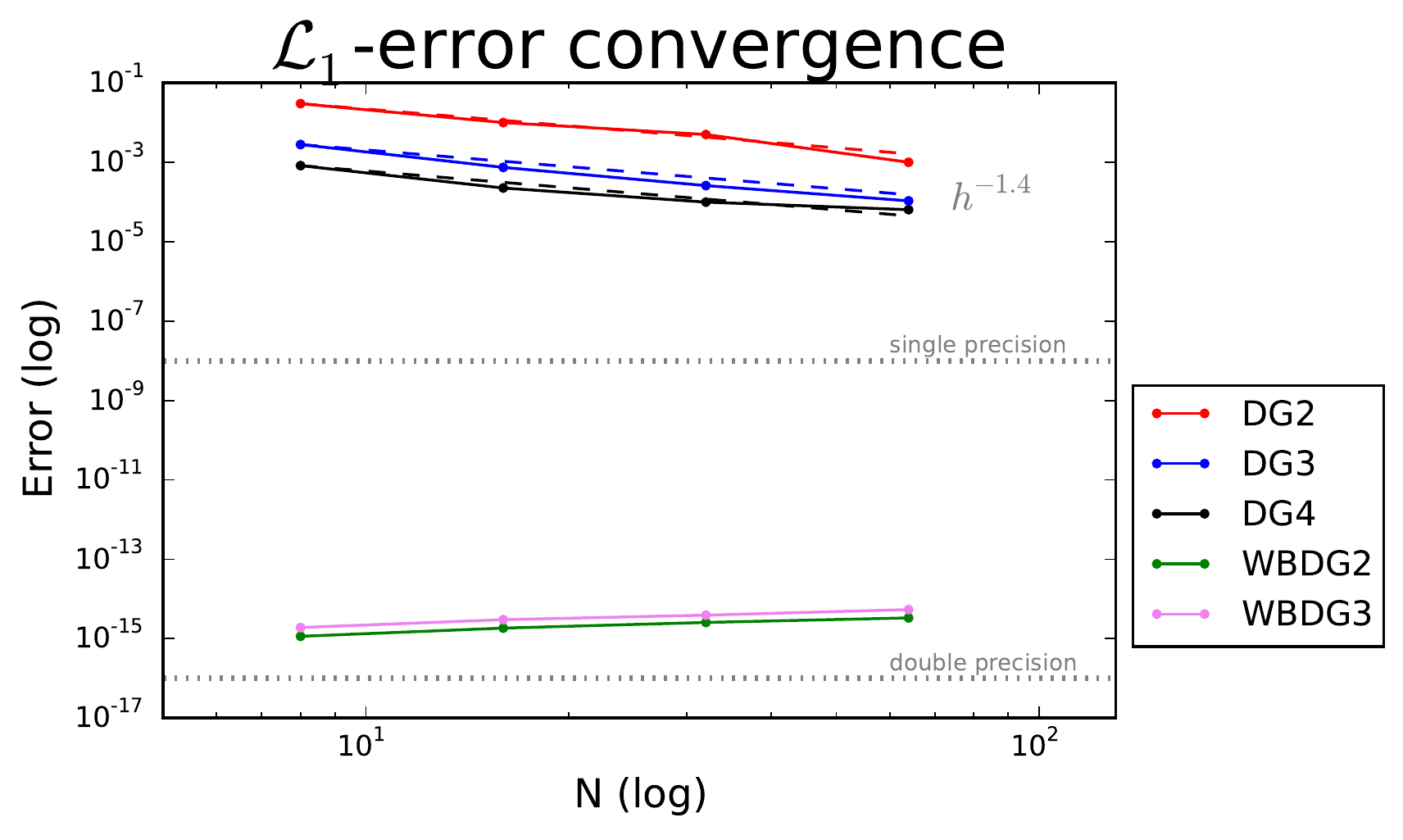}
         \caption{$\mathcal{L}_1$-error convergence for the 2-dimensional modified gresho vortex (see Eq.~\ref{eq:2dgresho}).}
         \label{fig:greshonorm}
\end{figure}

We consider a modified gresho vortex, where the pressure is modified to balance exactly a gravity source term. The initial conditions for the primitive variables are:
\begin{equation}
\label{eq:2dgresho}
\rho = 1.0, \quad v_x = -v_\theta \frac{(y-y_c)}{r}, \quad v_y = v_\theta \frac{(x-x_c)}{r}, \quad p = p(r),
\end{equation}
with the cross-radial velocity $v_\theta$ and pressure $p$:
\begin{equation*}
    v_{\theta}(r) =  \begin{cases}
    5r & r < 0.2 \\
    2-5r & 0.2\leq r < 0.4 \\
    0 & r\geq 0.4
  \end{cases}
\end{equation*}

\begin{equation*}
    p(r) =  \begin{cases}
    5 + \frac{25}{2}r^2 - \alpha \Phi & r < 0.2 \\
    9 - 4\log(0.2) + \frac{25}{2}r^2 - 20r + 4\log(r) - \alpha \Phi & 0.2\leq r < 0.4 \\
    3 + 4\log(2) - \alpha \Phi & r\geq 0.4
  \end{cases}
\end{equation*}

where $\alpha = 0.01$ and $\Phi = \frac{1}{r}$. One can easily verify that adding a gravity source term with the potential  $\alpha \Phi$ recovers the gresho vortex analytically and that this is a steady state solution of the Euler system. We run this initial condition until $T = 1.0$. In figure~\ref{fig:greshonorm} we can see the measured empirical $\mathcal{L}_1$-norm between the well balanced discretisation and the traditional discontinuous Galerkin method. We recover the expected convergence rate is of $\mathcal{O}(1.4)$ for the traditional Gresho vortex test case. 
Including the gravity source term does not alter the convergence properties of the scheme. We also see that the well-balance scheme maintain the equilibrium state to machine precision accuracy. 

\subparagraph{Simplified protoplanetary disc}
In the context of planet formation, it is customary to consider a stationary disc rotating around a single star, which is a steady state solution of the Euler-Poisson equations. 

In this paper, we consider a constant density disc defined in a $[-6,6]\times[-6,6]$ box, with the following initial conditions:
\begin{equation}
\begin{split}
\rho_{eq} = 1.0, \quad 
u_{eq} = -\frac{v_\theta}{r}y , \quad
v_{eq} = \frac{v_\theta}{r}x, \quad
p_{eq} = c_s^2\rho_{eq},
\end{split}
\end{equation}
where $v_\theta = \sqrt{\frac{1}{r}\left( 1 -\alpha^2\right)}$ is the orbital velocity (slightly sub-Keplerian), $c_s = \alpha v_K$ the speed of sound, given by the product of the Keplerian velocity $v_K$ and the disk aspect ratio $\alpha=0.03$, and the gravity potential of a unit point mass given by $\Phi = -\frac{1}{r}$.

We now describe in details how we set up our boundary conditions, for which great care is required in order to preserve the correct geometry of the problem and to stabilise the solution:

\begin{itemize}
\item{For the domain boundary conditions (on the box $[-6,6]\times[-6,6]$), the steady state solution is just imposed in ghost elements, as shown in \eqref{bcs}.}
\item{To minimise spurious effects due to the rotation of the disk near the end of the box domain, the constant density field $\rho_{eq}$ is multiplied with a tampering function $d(r)$. The following tampering function is taken:
\[d(r) = \frac{1}{1+\left(\frac{r}{r_0}\right)^{q}},\]
setting $q = 20$, $r_0 = 4.2$. This function was adopted after several other functions have been tried.  
Note that for the stability of the RKDG method it is important to consider functions which have well behaved derivatives at all orders. 
Another good candidate we have tried is the sigmoid function (not shown here).}
\item{The disc is an isolated system with no mass inflow and a (tampered) sharp edge. We need to introduce a buffer region near the disc edge where propagating waves are damped   
to reduce wave reflection. We use a methodology similar to \cite{valborro}, which smoothly relaxes the numerical solution to the equilibrium solution at the edge of the buffer zone using a function $R(r)$ so that
\[ \tilde{H}(u) = H(u) R(r) \]
Note that this function must leave the solution unaltered outside of the buffer region. 
In our experiments we set $R(r) =  \frac{1}{1+\exp(r^2-15.0)}$ where the parameter $15.0$ was chosen to set the size of the buffer region. 
In \cite{valborro}, the authors used a parabolic function $R(r)$ instead.}
\item{Similarly, at the centre of the disk, around $r = 0.0$, we use an inner buffer region where the numerical solution is set to the steady state solution. 
An inner radius of $r < 0.75$ is considered for the size of the inner buffer region.}
\end{itemize}
A perturbation is the added to the gravity field of the star. Physically, this perturbation can be interpreted as a planet. As such, the magnitude of the gravitational force exerted by the planet is very small in comparison to the gravitational force exerted by the star. We introduce this perturbation in the second term of equation \eqref{eq:planet2d}.
\begin{equation}
\label{eq:planet2d}
\begin{split}
    \nabla \Phi(\bm{x}) &=  \frac{\bm{x}}{(r^2+\epsilon^2)^{\frac{3}{2}}} + \eta\frac{\bm{x} - \bm{x}_p}{(r_p^2+\epsilon^2)^{\frac{3}{2}}}
    ,
\end{split}
\end{equation}
where  $r_p = \sqrt{||{\bm x}-\bm{x_p}||}$, $\bm{x}_p$ denotes the position of the perturbation

\[
\bm{x}_p =  \begin{pmatrix}
  x_p  \\
  y_p 
 \end{pmatrix} = \begin{pmatrix}
  r_c \cos \left(\frac{v_K}{r_c}t \right)  \\
  r_c \sin \left(\frac{v_K}{r_c}t \right) 
 \end{pmatrix},\]

fixed to be a circular orbit at $r=2.2$ with Keplerian velocity $v_K$ and $\epsilon = 0.01$ is the softening length for the planet. By varying $\eta$, we control the size of the perturbation. We test different sizes of $\eta$ to denote different sized planets, namely, $\eta = 3.1 \times 10^{-6},~9 \times 10^{-5}$ and $ 9.5 \times 10^{-4}$ which correspond to Earth, Neptune and Jupiter sized planets, respectively.

The system is evolved until 20 rotations are performed at $r = 2.2$, corresponding to approximately $T = 410$ in our normalised units. 

\begin{figure}
    \includegraphics[width=1.0\textwidth]{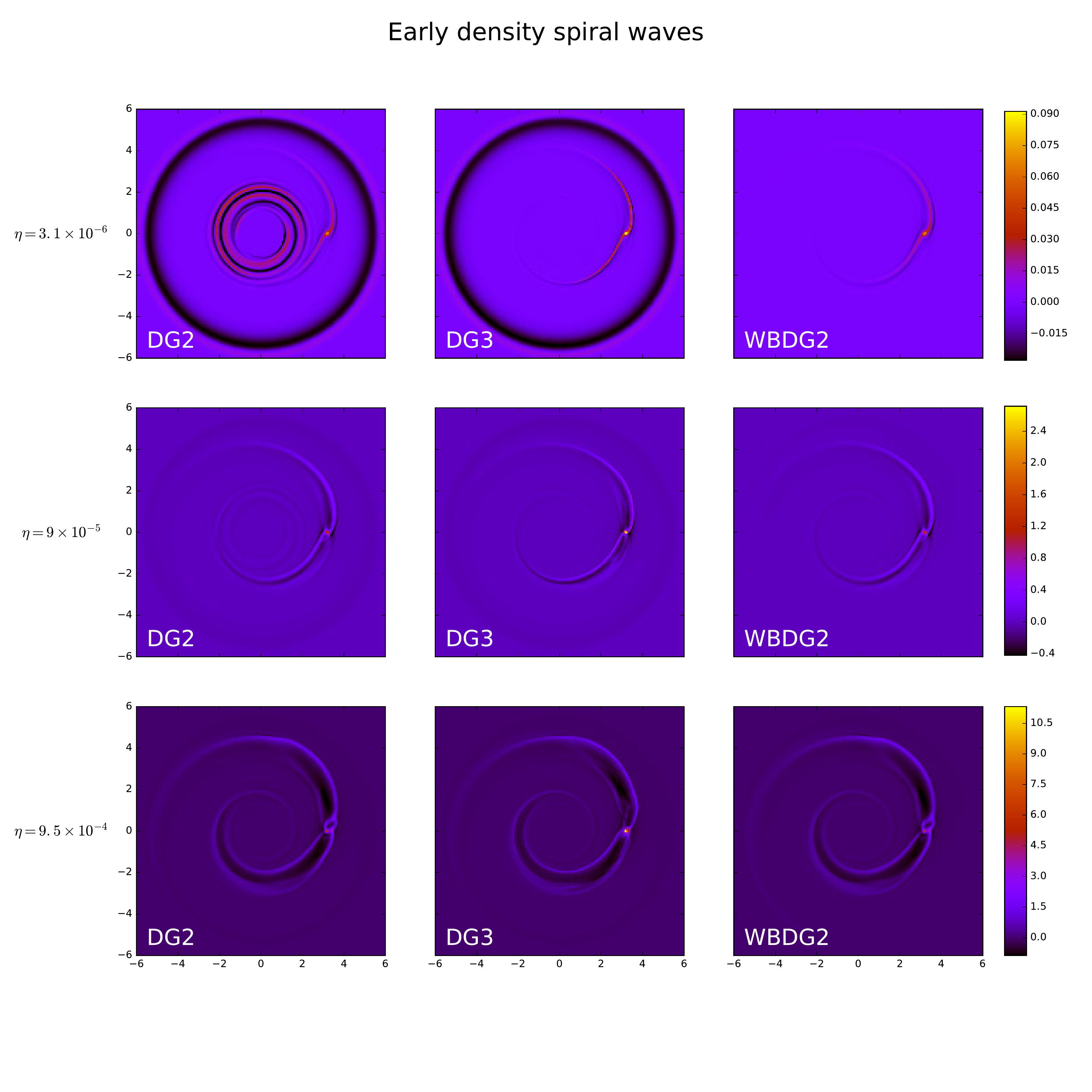}
    \caption{Density perturbations for non well-balanced method versus well-balanced method from dynamic equilibrium for varying perturbation of sizes on the gravity field, after 1 rotation, at approximately $T = 21$.}\label{fig:rot1}
\end{figure}

\begin{figure}
    \includegraphics[width=1.0\textwidth]{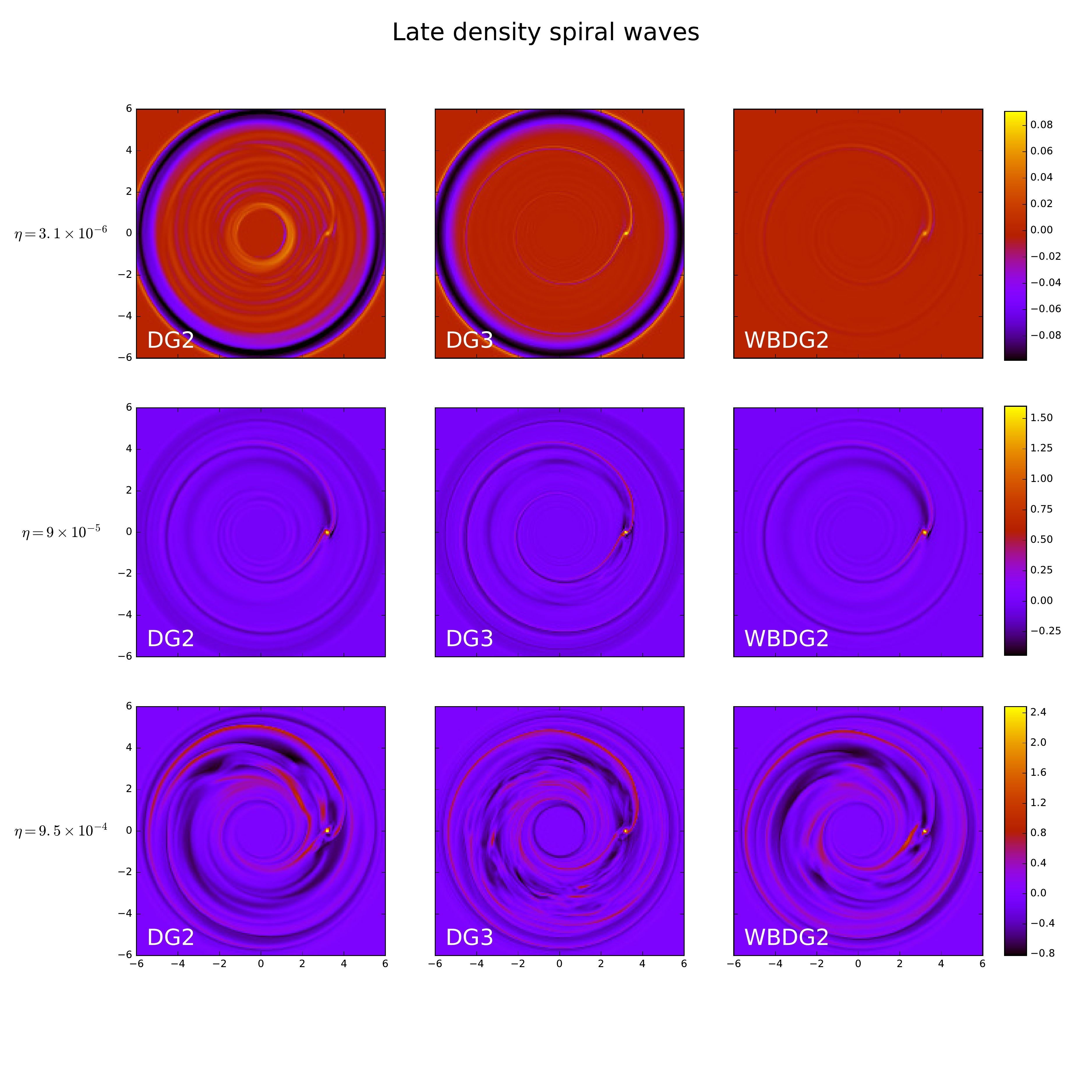}
    \caption{Density perturbations for non well-balanced method versus well-balanced method from dynamic equilibrium for varying perturbation of sizes on the gravity field, after 10 rotations, at approximately $T = 210$.}\label{fig:rot10}
\end{figure}

The results after the planet has performed only one rotation can be seen in figure~\ref{fig:rot1} and after 10 rotations in figure~\ref{fig:rot10}. For the smallest perturbation ($\eta = 3.1 \times 10^{-6}$), we note that already after one full rotation, the solution of the DG2 method has interacted with the waves generated by the mismatch between the inner boundary condition and the evolved solution, and this effect disappears when increasing the method to $3^{rd}$ order or when using the WBDG2 scheme. After 10 rotations it's clear that the perturbation has been lost in the numerical errors when using DG2, whereas for DG3 the solution remains very clean, both in the perturbation and the steady state background solution. Similarly, when using WBDG2, we observe a very clean perturbation on top of the unperturbed steady state background, although the resolution on the perturbation is lower than in the DG3 case. 

A similar behaviour is observed for the medium amplitude perturbation ($\eta = 9 \times 10^{-5}$), after one rotation. After 10 rotations, although the spiral density wave can be seen in all methods, both when using DG2 and DG3, artefacts are observed in the gap opened by the planet, whereas when using WBDG2 the gap remains cleaner. For the larger perturbation ($\eta = 9.5 \times 10^{-4}$), even though the effect from the boundary is still present, we see virtually no difference between DG2 and WBDG2. For such large sized planets, it is expected for a gap to be carved in the disc, and the regime of the study is very different. Indeed, after 10 rotations the disk is visibly unstable, and the solution has deviated enough from the steady state background solution that there's virtually no difference between DG2 and WBDG2. Indeed, for the simulation to reach 10 rotations, we had to stabilise all methods by using a positivity preserving limiter. Note that the large amplitude case is particularly interesting, because it demonstrates that our well-balanced scheme is robust enough to sustain large deviations from the adopted equilibrium state, recovering the properties of the corresponding non-well-balanced scheme.

Lastly, as denoted in table~\ref{tab:planetperts}, we show the time to solution required for different non well-balanced and well-balanced methods. In this example, it becomes clear the advantage of using a well-balanced scheme for long term evolution of small perturbations, as the necessary increase in order and resolution in the non well-balanced case can be translated into much longer simulation times.

\begin{table}[h]
\centering
\caption{Time to solution for the protoplanetary disc case after 10 rotations, for varying planet sizes.}
\label{tab:planetperts}
\begin{tabular}{||c | c c c ||} 
 \hline
 $\eta$ & DG2& DG3 & WBDG2 \\ [0.5ex] 
 \hline\hline
 $3.1 \times 10^{-6}$ & 1h42m & 16h30m6 & 3h18m \\ 
 $9 \times 10^{-5}$ & 1h42m & 16h10m & 3h21m \\ 
 $9.5 \times 10^{-4}$ & 1h30m & 15h30m & 3h30m \\  [1ex] 
 \hline
\end{tabular}
\end{table}

\section{Conclusion}
\label{sec:conclusion}

The motivation of this paper was to address the following three research questions:

\begin{itemize}
\item  \textbf{RQ 1:} Are there cases where using a high order scheme is sufficient to capture solutions close to a steady state?
\item  \textbf{RQ 2:} Under which circumstances is it necessary to use a well balanced method?
\item  \textbf{RQ 3:} What is the cost associated to each approach and how does it balance with accuracy?  
\end{itemize}

To address these questions, we compared a classical RKDG scheme with a novel well-balanced RKDG scheme using different numerical examples. We study the performance of these two numerical schemes under the regime of hydrostatic equilibrium and \textit{dynamic} equilibrium. This last type of equilibrium is of interest for many studies and simulations of proto-planetary systems. A summary of our results from the numerical experiments shown in section \ref{sec:results} follows:
\begin{itemize}
\item{When considering \textbf{hydrostatic equilibrium in one space dimension}, the non-well balanced high order method behaves very well. In particular, for waves with amplitude larger than the scheme's truncation error, the method was able to resolve the perturbation accurately as  expected. For example, when using the $3^{rd}$ or $4^{th}$ order method, we were able to reduce the truncation error down to $\mathcal{O}(10^{-8})$ and  $\mathcal{O}(10^{-12})$ for a resolution of $N = 64$, respectively.}
\item{Studying the \textbf{hydrostatic equilibrium case in a 2-dimensional setting}, we were able to reduce the error to $\mathcal{O}(10^{-8})$ and $\mathcal{O}(10^{-12})$ only for a resolution of $N_x = 64, N_y = 64$ when using a $3^{rd}$ or $4^{th}$ order method, respectively. However, it was observed that when using a well-balanced scheme, the required resolution (either in space or in polynomial degree) could be lowered without affecting the ability of the scheme to capture the waves.}
\item{Considering a \textbf{steady state with a non-trivial velocity in a 1-dimensional setting} (\emph{dynamic} equilibrium), we are able to reduce the error to $\mathcal{O}(10^{-6})$ and $\mathcal{O}(10^{-10})$ for a resolution of $N = 64$, for a $3^{rd}$ or $4^{th}$ order, respectively. When considering the same initial condition with a smaller perturbation in the pressure field, high order methods appear to fail in capturing properly the wave (in particular for perturbation $\eta = 1\times10^{-8}$). Only our well-balanced scheme is robust enough to capture the wave dynamics properly.}
\item{For the \textbf{modified Gresho vortex, a 2-dimensional steady state solution with a non-trivial velocity}, we observe that the non well-balanced method converges with the expected numerical accuracy, even in the presence of the source terms. The well-balanced scheme is preserving the stationary solution down to machine precision accuracy.}
\item{For the \textbf{idealised disc-planet interaction case},  a 2-dimensional steady state solution with a non-trivial velocity, we observe that for small and medium amplitude perturbations, the well-balanced scheme is clearly superior, being able to keep the background solution for longer than the non well-balanced methods, although the increase in the approximation order improves the solution significantly. In particular, for the medium amplitude perturbation we observe the opening of a shallow gap. In both DG2 and DG3 this gap is not very clean. In the context of planet migration, the gap opening is important because it decreases the angular momentum exchange between the planet and disc, which in turn is translated into a decrease of the migration rate of the planet \cite{gapopen}. For both the small and medium amplitude perturbations, we are able to observe the expected density spiral arms without the strong numerical artefacts that we see when using classical  high order schemes. Indeed, in this case the density wave seems to be corrupted mostly from the inner and outer boundaries. In order to stabilise the solution when using DG2 and DG3, in addition to sophisticated boundary conditions, we have to use a positivity preserving limiter \cite{pp}. 
For the large amplitude perturbation, we see no difference between DG2 and WBDG2. Our hypothesis is that the solution can't be represented as a simple superposition of a steady state solution and a time dependent small perturbation, however, we note that the WBDG2 scheme behaves as DG2 (also requiring the positivity preserving limiter), which points towards our well-balanced method being robust for large perturbations.}
\end{itemize}

The use of a well-balanced scheme for the rotating disk case seems a good compromise (\textbf{RQ2}), while for simpler 1-dimensional hydrostatic equilibrium problems, one can beat down the truncation error fairly easily by raising either the order or the resolution of the scheme, without using the well balance correction (\textbf{RQ1}).

We note that the well-balanced correction does not come without a cost (\textbf{RQ3}). As shown in Figure~\ref{fig:timing_nomem}, the well-balanced correction can slow down the code significantly. This can be alleviated by pre-computing and storing all the variables from the steady state solution. However, this means that the memory requirements for this algorithm almost double (in comparison to the classical RKDG scheme). Due to the compute intensity of DG methods, GPUs are usually the appropriate hardware to run these methods  \cite{kridovonova}, which are often limited in memory so this is something worth considering when choosing the appropriate implementation. 

Note than one very restrictive condition for our well-balanced scheme is to know the exact form of the equilibrium solution everywhere, either in analytical or tabulated form. 
Moreover, a high-order scheme, if it is robust enough to capture the wave dynamics, will always deliver higher accuracy than a well-balanced, low order scheme. 
In \cite{david2018}, we show that the DG method can be a competitive method in planet-disc interaction studies, even without the well balanced correction, if one uses enough grid points, a conservative slope limiter and a careful boundary condition strategy.

In conclusion, we have shown that the well balanced property in numerical schemes is important for two reasons: 1- we are able to capture low amplitude waves propagating in non-trivial equilibrium states without resorting to complex boundary conditions strategies and 2-  we are able to solve for very small perturbations using lower order methods, which requires significantly less computations per time step, in particular when considering multi-dimensional problems. Both these points are relevant for setups like the one discussed here, namely the long term evolution of the sightly perturbed multi-dimensional equilibrium disc solution. As further steps, we hope to test and potentially study an extension of this scheme to capture equilibrium solution with discontinuities and to consider arbitrary general equilibrium states along the lines of \cite{Ricchiuto15}.

\section*{Acknowledgments}
MHV is supported by the UZH Candoc Scholarship. The computing resources were provided by the S3IT cluster at University of Zurich. The authors would like to thank the anonymous reviewers for their valuable comments and suggestions to improve the quality of the paper.

\appendix
\section{One dimensional moving steady state solution}
\label{Appendix1}
In this section we are interested in the construction of a simple 1-dimensional test case which is a steady state solution to the 1-dimensional Euler equations with a non trivial velocity field. 
The objective is to find the quartet of functions $(\rho, v, p, \phi)$ such that they fulfil the following:

\begin{subequations}
\label{eq:1deuler}
\begin{align}
&\frac{\partial}{\partial x} (\rho v) = 0 \label{eq:euler1d1} \\
&\frac{\partial}{\partial x} (\rho v^2 + p) = -\rho \frac{\partial}{\partial x} \Phi \label{eq:euler1d2} \\
&\frac{\partial}{\partial x}\big((E+p)v\big) = -\rho v \frac{\partial}{\partial x} \Phi \label{eq:euler1d3}
\end{align}
\end{subequations}

From \eqref{eq:euler1d1}, we have $\rho v = \mbox{const}$, whereas for \eqref{eq:euler1d2} and \eqref{eq:euler1d3}:
\begin{align*}
\rho v\frac{\partial}{\partial x}  v + \frac{\partial}{\partial x} p &= -\rho \frac{\partial}{\partial x} \Phi \\
\frac{\partial}{\partial x}  \big((E + p)v\big) &= -\rho v \frac{\partial}{\partial x} \Phi
\end{align*}

Noting that $E = \frac{p}{\gamma -1} + \frac{1}{2}\rho v^2$, \eqref{eq:euler1d3} yields:

\[\rho v^2 \frac{\partial}{\partial x} v + \frac{\partial}{\partial x}\left(p v \frac{\gamma}{\gamma -1}\right) = -\rho v \frac{\partial}{\partial x} \Phi
\]

Substituting \eqref{eq:euler1d2} into \eqref{eq:euler1d3}, one can solve find $p$ if we assume some form for $\rho$ (and consequently for $v$).

Setting $\rho = \exp(-x)$, thus $v=\exp(x)$ and $p = exp(-\gamma x)$. An expression for $\Phi$ can be written by solving the differential equation in \eqref{eq:euler1d2}, yielding:
$\frac{\partial}{\partial x}\Phi = \exp(x)(-\exp(x)+\gamma\exp(-\gamma x))$.

\section{A simple equilibrium solution for proto-planetary discs}
\label{Appendix2}

The orbital speed for a gas can be calculated from the Euler-Poisson equations:
\[ \frac{\partial \textbf{v}}{\partial t} + (\textbf{v}\cdot \nabla)\textbf{v}  = -\frac{1}{\rho} \nabla p - \nabla \Phi , \]
where $p$ denotes the pressure, $\rho$ the density and $\Phi$ the gravitational potential. 
We can rewrite the second term as \footnote{Using the following identity $\nabla (\textbf{A}\cdot \textbf{B}) = \textbf{A}\times (\nabla \times \textbf{B}) - (\nabla \times \textbf{A}) \times \textbf{B} + (\textbf{A}\cdot \nabla) \textbf{B} + (\textbf{B}\cdot \nabla) \textbf{A}$}:

\[ (\textbf{v} \cdot\nabla) \textbf{v} = \frac{1}{2}\nabla \textbf{v}^T \textbf{v} - \textbf{v} \times(\nabla \times \textbf{v} )\]

Assuming a steady state ($\frac{\partial}{\partial t} = 0$) axisymmetric solution, we derive the orbital velocity:

\[ \frac{v_{\phi}^2}{r} = \frac{1}{\rho}\frac{\partial }{\partial r} p + \frac{\partial}{\partial r} \Phi \]

Furthermore, defining the Keplerian velocity: $v_{\rm K}(r) = \sqrt{ r \frac{\partial}{\partial r} \Phi}$ and the constant disk aspect ratio \cite{armitageplanet} with
\[\alpha = \frac{\sqrt{\frac{p(r)}{\rho(r)}}}{v_{\rm K}(r)} \, ,\]
we can deduce the relation for the pressure to be $p(r) = \alpha^2 \rho(r)v_{\rm K}^2$. Finally, we obtain the equilibrium tangential velocity $v_\phi$ knowing the constant $\alpha$ and the profile $\rho(r)$.

\section{Supplementary results}
\label{ap:extraresults}

In this section we provide the time to solution for the numerical experiments performed in section~\ref{sec:results} and convergence plots associated to the perturbation tests.

\subsection{1-dimensional hydrostatic}
For convenience, we restate the initial conditions: an ideal gas $\gamma = 1.4$ in isothermal equilibrium state and a linear gravitational potential $\Phi_x = gx$ is considered:

\begin{equation}
\begin{split}
\label{eq:1dhydrostatic_a}
    &\rho_{eq}(x) = \rho_0 \exp \bigg(-\frac{\rho_0 g}{p_0}x \bigg) \\
    &u_{eq}(x) = 0\\
    &p(x,t=0) = p_{eq}(x) + \eta \exp \bigg(-\frac{\rho_0 g}{p_0}\frac{(x-0.5)^2}{0.01}\bigg) \\
\end{split}
\end{equation}
with $\rho_0 = 1.0$, $p_0 = 1.0$ and $g=1.0$.

\begin{table}[h]
\centering
\caption{Time to solution for the 1-dimensional hydrostatic equilibrium \eqref{eq:1dhydrostatic_a} (s) for perturbation size $\eta = 1\times10^{-2}$, $1\times10^{-4}$ and $1\times10^{-6}$, respectively. }
\label{tab:hydrostatic1-6}

\begin{tabular}{||c | c c c c ||} 
 \hline
 N & DG1& DG2 & WBDG2 & WBDG3 \\ [0.5ex] 
 \hline\hline
 8 & 0.37 & 0.56 & 0.05 & 0.14\\ 
 16 & 0.62 & 1.28 & 0.10 & 0.29\\ 
 32 & 2.05 & 4.79 & 0.20 & 0.69\\ 
 64 & 12.8 & 32.4 & 0.62 & 3.25\\
 128 & 103 & 270 & 3.84 & 24.1\\   [1ex] 
 \hline
\end{tabular} \begin{tabular}{||c | c c c c ||} 
 \hline
 N & DG2& DG3 & WBDG2 & WBDG3 \\ [0.5ex] 
 \hline\hline
 8 & 0.37 & 0.56 & 0.05 & 0.14\\ 
 16 & 0.62 & 1.28 & 0.10 & 0.29\\ 
 32 & 0.18 & 4.79 & 0.20 & 0.69\\ 
 64 & 12.8 & 32.4 & 0.62 & 3.25\\
 128 & 103 & 270 & 3.84 & 24.1\\   [1ex] 
 \hline
\end{tabular}\\   [1ex] 
\begin{tabular}{||c | c c c c ||} 
 \hline
 N & DG2& DG3 & WBDG2 & WBDG3 \\ [0.5ex] 
 \hline\hline
 8 & 0.37 & 0.56 & 0.05 & 0.14\\ 
 16 & 0.62 & 1.28 & 0.10 & 0.29\\ 
 32 & 2.05 & 4.79 & 0.20 & 0.69\\ 
 64 & 12.8 & 32.4 & 0.62 & 3.25\\
 128 & 103 & 270 & 3.84 & 24.1\\   [1ex] 
 \hline
\end{tabular}
\end{table}

\subsection{2-dimensional hydrostatic}

\begin{figure}
    \centering
    \begin{subfigure}[b]{0.49\textwidth}
        \includegraphics[width=\textwidth]{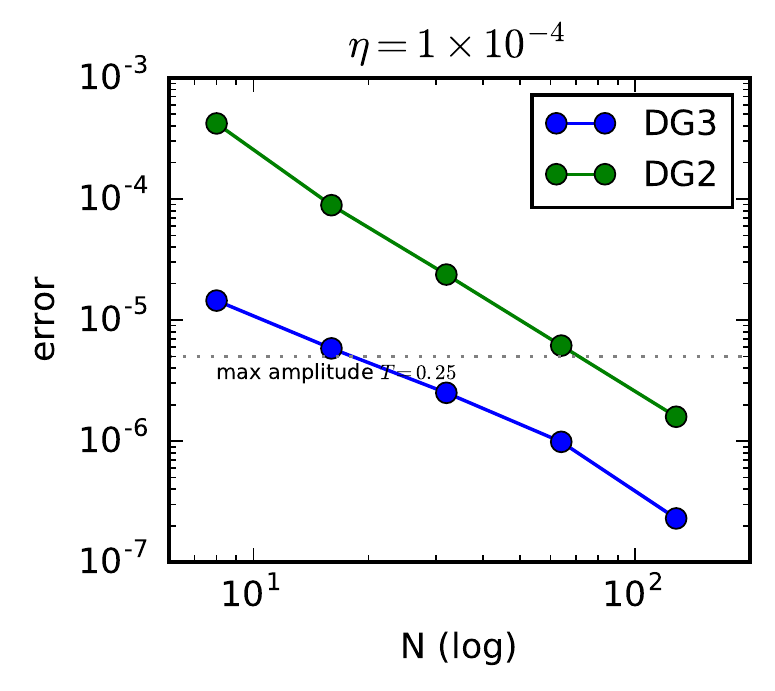}
            \end{subfigure}
                \begin{subfigure}[b]{0.49\textwidth}
        \includegraphics[width=\textwidth]{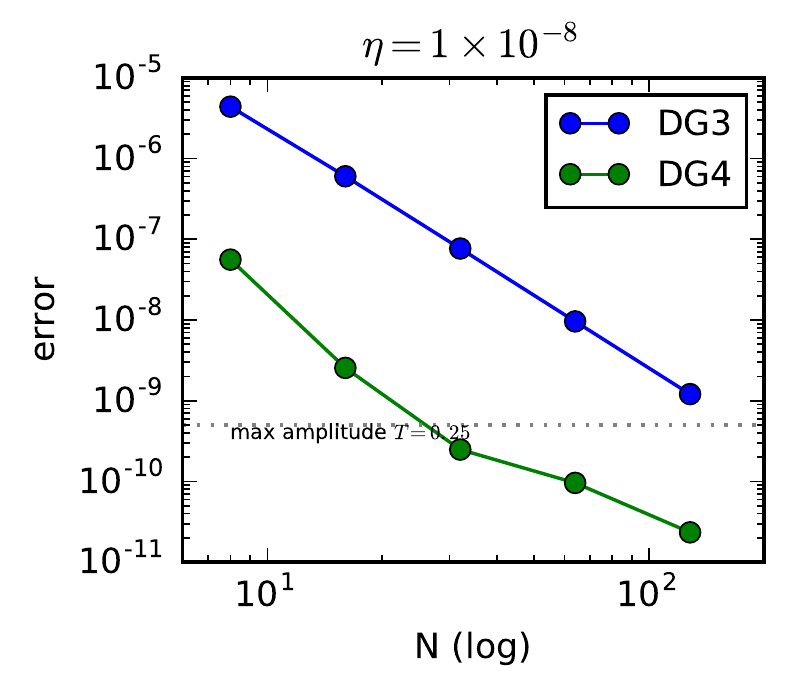}
    \end{subfigure}
     \caption{Non well-balanced method vs well-balanced method for hydrostatic equilibrium with varying amplitude perturbation on the pressure field for initial conditions \eqref{eq:2dhydrostatic_a}.}
    \label{fig:conv2dpulse_a}
\end{figure}

Ideal gas $\gamma = 1.4$, in isothermal equilibrium and a linear gravitational potential $\Phi = g(x + y)$. Unit square domain $\bm{x}\in[0,1]\times[0,1]$:
\begin{equation}
\label{eq:2dhydrostatic_a}
\begin{split}
    \rho_{eq}(x,y) &= \rho_0 \exp \bigg(-\frac{\rho_0 g}{p_0}(x+y) \bigg) \\
    u_{eq}(x,y) &= 0\\
    v_{eq}(x,y) &= 0\\
    p_{eq}(x,y) &= p_0 \exp \bigg(-\frac{\rho_0 g}{p_0}(x+y) \bigg) \\
\end{split}
\end{equation}
with $\rho_0 = 1$, $p_0 = 1$ and $g=1$. The time to solution is shown on tables~\ref{tab:hydrostatic2dtim},~\ref{tab:2dperts_a} and error convergence plots in figure~\ref{fig:conv2dpulse_a}.

\begin{table}[h]
\centering
\caption{Time to solution for hydrostatic equilibrium \eqref{eq:2dhydrostatic_a} (s) at $T = 10.0$.}
\label{tab:hydrostatic2dtim}
\begin{tabular}{||c | c c c c c ||} 
 \hline
 $N_x$ & DG2& DG3 & DG4 & WBDG2 & WBDG3 \\ [0.5ex] 
 \hline\hline
 8 & 1.53 & 3.36 & 10.8 &1.85 & 4.63 \\ 
 16 & 2.92 & 7.32 & 24.6 & 3.74 & 9.52\\ 
 32 & 6.88 & 19.9 & 82.1 & 8.74 & 24.7\\ 
 64 & 19.3 & 101 & 525 & 2.46 & 130\\
 128 & 119 & 777 & 4310 & 154 & 981\\ [1ex] 
 \hline
\end{tabular}
\end{table}

\begin{table}[h]
\centering
\caption{Time to solution for hydrostatic equilibrium  \eqref{eq:2dhydrostatic_a} (s) for perturbation sizes $ \eta = 1\times10^{-4}$, $1\times10^{-8}$ at $T = 0.25$.}
\label{tab:2dperts_a}
\begin{tabular}{||c | c c c c ||} 
 \hline
 $N_x$ & DG2& DG3 & WBDG2 & WBDG3 \\ [0.5ex] 
 \hline\hline
 8 & 0.04 & 0.08 & 0.05 & 0.10\\ 
 16 & 0.08 & 0.17 & 0.09 & 0.22\\ 
 32 & 0.17 & 0.49 & 0.22 & 0.64\\ 
 64 & 0.50 & 2.51 & 0.69 & 3.26\\
 128 & 2.98 & 19.5 & 3.84 & 24.5\\   [1ex] 
 \hline
\end{tabular}
\begin{tabular}{||c | c c c c ||} 
 \hline
 $N_x$ & DG3& DG4 & WBDG2 & WBDG3 \\ [0.5ex] 
 \hline\hline
 8 & 0.08 & 0.27 & 0.03 & 0.10\\ 
 16 & 0.17 & 0.62 & 0.07 & 0.22\\ 
 32 & 0.50 & 2.06 & 0.17 & 0.62\\ 
 64 & 2.50 & 13.2 & 0.61 & 3.26\\
 128 & 19.5 & 108 & 3.86 & 24.5\\   [1ex] 
 \hline
\end{tabular}
\end{table}

\subsection{1-dimensional dynamic}
\begin{figure}
    \centering
    \begin{subfigure}[b]{0.49\textwidth}
        \includegraphics[width=\textwidth]{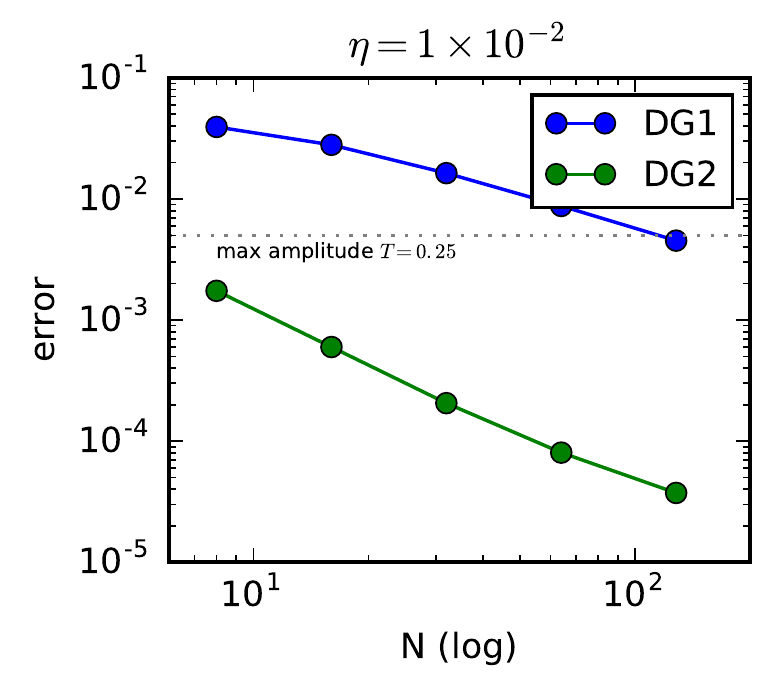}
            \end{subfigure}
                \begin{subfigure}[b]{0.49\textwidth}
        \includegraphics[width=\textwidth]{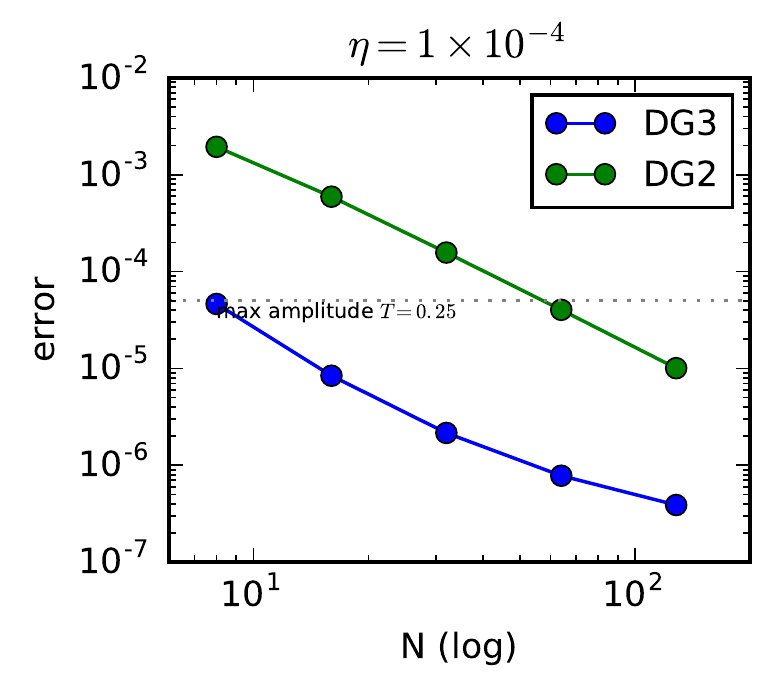}
    \end{subfigure}\\
    \begin{subfigure}[b]{0.49\textwidth}
        \includegraphics[width=\textwidth]{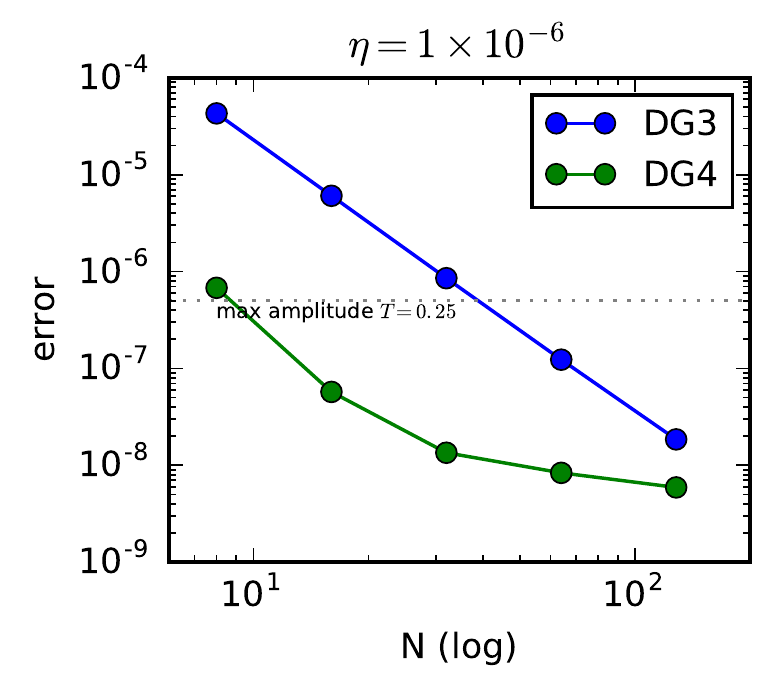}
    \end{subfigure}
    \begin{subfigure}[b]{0.49\textwidth}
        \includegraphics[width=\textwidth]{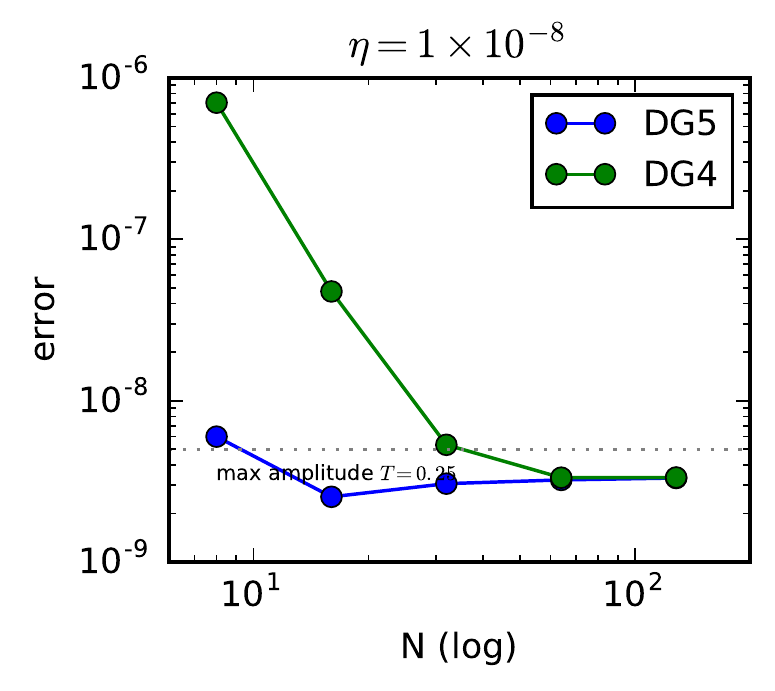}
    \end{subfigure}
    \caption{Non well-balanced method versus well-balanced method for hydrostatic equilibrium with varying amplitude perturbation on the pressure field for initial conditions \eqref{eq:1ddynamic_a}.}
    \label{fig:conv1ddympulse}
\end{figure}

Ideal steady gas $\gamma = 1.4$ with a nonzero velocity field and a non linear gravitational field. 
\begin{equation}
\label{eq:1ddynamic_a}
\begin{split}
    \rho_{eq}(x) &= \rho_0 \exp \bigg(-\frac{\rho_0 g}{p_0}x \bigg) \\
    u_{eq}(x) &= \exp(x)\\
    p_{eq}(x) &= \exp \bigg(-\frac{\rho_0 g}{p_0}x \bigg)^{\gamma} \\
\end{split}
\end{equation}
with $\rho_0 = 1$, $p_0 = 1$ and a non linear potential $\phi = \exp(x) (-\exp(x) + \gamma \exp(-\gamma x))$. The time to solution for this test case is very similar to the 1-dimensional hydrostatic equilibrium case, and is thus omitted. The error convergence plots are shown in figure~\ref{fig:conv1ddympulse}.

\clearpage
\bibliographystyle{plain}
\bibliography{template}

\begin{thebibliography}{10}

\bibitem{ern}
Daniele Antonio Di~Pietro and Alexandre Ern.
\newblock {\em Mathematical Aspects of Discontinuous Galerkin Methods},
  volume~69.
\newblock 01 2012.

\bibitem{armitageplanet}
Philip~J. Armitage.
\newblock Dynamics of protoplanetary disks.
\newblock {\em Annual Review of Astronomy and Astrophysics}, 49(1):195--236,
  2011.

\bibitem{bermudez2017}
Alfredo Berm\'{u}dez, Xi\'{a}n L\'{o}pez, and M.~Elena V\'{a}zquez-Cend\'{o}n.
\newblock Numerical solution of non-isothermal non-adiabatic flow of real gases
  in pipelines.
\newblock {\em J. Comput. Phys.}, 323(C):126--148, October 2016.

\bibitem{bermudez94}
Alfredo Berm\'udez and Maria~Elena Vazquez.
\newblock Upwind methods for hyperbolic conservation laws with source terms.
\newblock {\em Computers \& Fluids}, 23(8):1049 -- 1071, 1994.

\bibitem{berthon2012}
Christophe Berthon and Fran\c{c}oise Foucher.
\newblock Efficient well-balanced hydrostatic upwind schemes for shallow-water
  equations.
\newblock 231:4993–5015, 06 2012.

\bibitem{pares2007}
Manuel~J. Castro, Alberto Pardo~Milan\'es, and Carlos Par\'es.
\newblock Well-balanced numerical schemes based on a generalized hydrostatic
  reconstruction technique.
\newblock {\em Mathematical Models and Methods in Applied Sciences},
  17(12):2055--2113, 2007.

\bibitem{Klingenberg}
Praveen Chandrashekar and Christian Klingenberg.
\newblock A second order well-balanced finite volume scheme for euler equations
  with gravity.
\newblock {\em SIAM Journal on Scientific Computing}, 37(3):B382--B402, 2015.

\bibitem{cs5}
Bernardo Cockburn and Chi-Wang Shu.
\newblock The runge-kutta discontinuous galerkin method for conservation laws
  v.
\newblock {\em J. Comput. Phys.}, 141(2):199--224, April 1998.

\bibitem{valborro}
M.~De~Val-Borro, R.~G. Edgar, P.~Artymowicz, P.~Ciecielag, P.~Cresswell,
  G.~D'Angelo, E.~J. Delgado-Donate, G.~Dirksen, S.~Fromang, A.~Gawryszczak,
  H.~Klahr, W.~Kley, W.~Lyra, F.~Masset, G.~Mellema, R.~P. Nelson, S.-J.
  Paardekooper, A.~Peplinski, A.~Pierens, T.~Plewa, K.~Rice, C.~Schäfer, and
  R.~Speith.
\newblock A comparative study of disc–planet interaction.
\newblock {\em Monthly Notices of the Royal Astronomical Society},
  370(2):529--558, 2006.

\bibitem{dedner2001}
A.~Dedner, I.L Sofronov, and M~Wesenberg.
\newblock Transparent boundary conditions for mhd simulations in stratified
  atmospheres.
\newblock {\em Journal of Computational Physics}, 171(2):448 -- 478, 2001.

\bibitem{kridovonova}
Martin Fuhry, Andrew Giuliani, and Lilia Krivodonova.
\newblock Discontinuous galerkin methods on graphics processing units for
  nonlinear hyperbolic conservation laws.
\newblock {\em CoRR}, abs/1601.07944, 2016.

\bibitem{gottlieb}
Sigal Gottlieb and Chi-Wang Shu.
\newblock Total variation diminishing runge-kutta schemes.
\newblock {\em Math. Comput.}, 67(221):73--85, January 1998.

\bibitem{Leroux96}
J.~M. Greenberg and A.~Y. Leroux.
\newblock A well-balanced scheme for the numerical processing of source terms
  in hyperbolic equations.
\newblock {\em SIAM Journal on Numerical Analysis}, 33(1):1--16, 1996.

\bibitem{pp}
Xiangyu~Y. Hu, Nikolaus~A. Adams, and Chi-Wang Shu.
\newblock Positivity-preserving method for high-order conservative schemes
  solving compressible euler equations.
\newblock {\em Journal of Computational Physics}, 242(Supplement C):169 -- 180,
  2013.

\bibitem{Kapelli}
{K\"appeli, R.} and {Mishra, S.}
\newblock A well-balanced finite volume scheme for the euler equations with
  gravitation - the exact preservation of hydrostatic equilibrium with
  arbitrary entropy stratification.
\newblock {\em A\&A}, 587:A94, 2016.

\bibitem{leveque2006}
R.J. LeVeque, O.~Steiner, A.~Gautschy, D.~Mihalas, E.A. Dorfi, and
  E.~M{\"u}ller.
\newblock {\em Computational Methods for Astrophysical Fluid Flow: Saas-Fee
  Advanced Course 27. Lecture Notes 1997 Swiss Society for Astrophysics and
  Astronomy}.
\newblock Saas-Fee Advanced Course. Springer Berlin Heidelberg, 2006.

\bibitem{solarbook}
L.A. McFadden, T.~Johnson, and P.~Weissman.
\newblock {\em Encyclopedia of the Solar System}.
\newblock Encyclopedia of the Solar System Series. Elsevier Science, 2006.

\bibitem{galaxybook}
H.~Mo, F.~van~den Bosch, and S.~White.
\newblock {\em Galaxy Formation and Evolution}.
\newblock Galaxy Formation and Evolution. Cambridge University Press, 2010.

\bibitem{Noelle}
Sebastian {Noelle}, Yulong {Xing}, and Chi-Wang {Shu}.
\newblock {High-order well-balanced finite volume WENO schemes for shallow
  water equation with moving water.}
\newblock {\em {J. Comput. Phys.}}, 226(1):29--58, 2007.

\bibitem{noelle2016}
Sebastian {Noelle}, Yulong {Xing}, and Chi-Wang {Shu}.
\newblock {High-order well-balanced schemes}.
\newblock In {\em {Numerical methods for balance laws}}, pages 1--66. Caserta:
  Dipartimento di Matematica, Seconda Universit\`a di Napoli, 2009.

\bibitem{Pares2006}
Carlos Par\'es.
\newblock Numerical methods for nonconservative hyperbolic systems: a
  theoretical framework.
\newblock {\em SIAM Journal on Numerical Analysis}, 44(1):300--321, 2006.

\bibitem{gapopen}
R.~R. Rafikov.
\newblock Planet migration and gap formation by tidally induced shocks.
\newblock {\em The Astrophysical Journal}, 572(1):566, 2002.

\bibitem{Ricchiuto15}
Mario Ricchiuto.
\newblock An explicit residual based approach for shallow water flows.
\newblock {\em J. Comput. Phys.}, 280(C):306--344, January 2015.

\bibitem{shu89}
Chi-Wang Shu and Stanley Osher.
\newblock Efficient implementation of essentially non-oscillatory
  shock-capturing schemes, ii.
\newblock {\em J. Comput. Phys.}, 83(1):32--78, July 1989.

\bibitem{surville2016}
Cl\'ement Surville, Lucio Mayer, and Douglas N.~C. Lin.
\newblock Dust capture and long-lived density enhancements triggered by
  vortices in 2d protoplanetary disks.
\newblock {\em The Astrophysical Journal}, 831(1):82, 2016.

\bibitem{supernova1996}
F.-K. {Thielemann}, K.~{Nomoto}, and M.-A. {Hashimoto}.
\newblock {Core-Collapse Supernovae and Their Ejecta}.
\newblock {\em Applied Physics Journal}, 460:408, March 1996.

\bibitem{david2018}
David~A. Velasco, Maria~Han Veiga, Fr\'ed\'eric Masset, and Romain Teyssier.
\newblock Planet-disc interactions with discontinuous galerkin methods using
  gpus.
\newblock {\em MNRAS (submitted)}, 42(2):641--666, 2004.

\bibitem{yee}
Wei {Wang}, Chi-Wang {Shu}, H.C. {Yee}, and Bj\"orn {Sj\"ogreen}.
\newblock {High order finite difference methods with subcell resolution for
  advection equations with stiff source terms.}
\newblock {\em {J. Comput. Phys.}}, 231(1):190--214, 2012.

\bibitem{Xing2013}
Yulong Xing and Xiangxiong Zhang.
\newblock Positivity-preserving well-balanced discontinuous galerkin methods
  for the shallow water equations on unstructured triangular meshes.
\newblock {\em Journal of Scientific Computing}, 57(1):19--41, Oct 2013.

\bibitem{shuconvergence}
Qiang Zhang and Chi-Wang Shu.
\newblock Error estimates to smooth solutions of runge--kutta discontinuous
  galerkin methods for scalar conservation laws.
\newblock {\em SIAM Journal on Numerical Analysis}, 42(2):641--666, 2004.

\bibitem{zhang}
Xiangxiong Zhang and Chi-Wang Shu.
\newblock Positivity-preserving high order discontinuous galerkin schemes for
  compressible euler equations with source terms.
\newblock {\em Journal of Computational Physics}, 230(4):1238 -- 1248, 2011.

\end{thebibliography}

\end{document}